\documentclass[12pt]{article}
    \usepackage{amsmath}
    \usepackage{amssymb}
    \usepackage{amsthm}

\usepackage{pstricks}

\usepackage[T2A]{fontenc}
\usepackage[cp1251]{inputenc}

\usepackage{graphicx,psfrag}

      \newcommand {\al}   {\alpha}          \newcommand {\bt}  {\beta}
      \newcommand {\gam } {\gamma}          \newcommand {\Gam} {\Gamma}
      \newcommand {\del}  {\delta}          
              \newcommand {\ve}  {\varepsilon}
                 \newcommand {\vphi} {\varphi}
      \newcommand {\lam}  {\lambda}         
      
                \newcommand {\Om}  {\Omega}
      \newcommand {\pl}   {\partial}        \newcommand {\s}   {\sigma}

      \newcommand {\RRR}  {{\mathbb R}}

      \newcommand {\ZZZ}  {{\mathbb Z}}     
\newcommand {\AAA}{\bar{A}}  \newcommand {\BBB}{\bar{B}} \newcommand {\OOO}{\bar{O}}
              
      \newcommand {\FFF}  {{\cal F}}        %\newcommand {\BBB}  {{\cal B}}
      \newcommand {\TTT}  {{\cal T}}

            \newcommand {\interval}{[-\pi/2,\, \pi/2]}
     \newcommand {\beq}  {\begin{equation}}
      \newcommand {\eeq}  {\end{equation}}
     \newcommand {\beqo}  {\begin{equation*}}
      \newcommand {\eeqo}  {\end{equation*}}

      \newcommand {\doubpar}  {\rput(0.08,-0.01){
      \scalebox{0.25}{\psecurve[linewidth=2pt]
      (-0.21875,-0.25)(-0.25,0)(-0.21875,0.25)(-0.125,0.5)(0,0.7071)(0.25,1)
      \psecurve[linewidth=2pt](0.21875,-0.25)(0.25,0)(0.21875,0.25)(0.125,0.5)(0,0.7071)(-0.25,1)
      \psline[linewidth=2pt](-0.25,0)(0.25,0)}}\hspace*{0.95mm}}

      \newcommand {\rr} {retroreflector}
      \newcommand {\rrs} {retroreflectors}
      \newcommand {\prr} {perfect retroreflector}
      \newcommand {\prrs} {perfect retroreflectors}
      \newcommand {\aprr} {asymptotically perfect retroreflector}
      \newcommand {\aprrs} {asymptotically perfect retroreflectors}
      \newcommand {\dopa} {helmet}
      \newcommand {\Dopa} {Helmet}

      \newtheorem{theorem}{Theorem}
      
      \newtheorem{utv}{Statement}
      \newtheorem{remark}{Remark}
      \newtheorem{opr}{Definition}

\author{Alexander Plakhov}

\title{Mathematical retroreflectors}

\date{University of Aveiro, Portugal}

\begin{document}

\maketitle

\begin{abstract}
Retroreflectors are optical devices that reverse the direction of incident beams of light. Here we present a collection of billiard type \rrs\ consisting of four objects; three of them are \aprrs, and the fourth one is a \rr\ which is very close to perfect. Three objects of the collection have recently been discovered and published or submitted for publication. The fourth object --- {\it notched angle} --- is a new one; a proof of its retroreflectivity is given.
\end{abstract}

\begin{quote}
\end{quote}

\begin{quote}
\end{quote}

\section{Introduction}\label{sect introduction}

%\subsection{Definition of a \rr}\label{subsec def retroreflector}

In everyday life, optical devices that reverse the direction of all (or a significant part of) incident beams of light are called {\it retroreflectors}. They are widely used, for example, in road safety. Some artificial satellites in Earth orbit also carry \rrs. We are mostly interested here in {\it \prrs} that reverse the direction of {\it any} incident beam of light to {\it exactly opposite}. An example of \prr\ based on {\it light refraction} is the Eaton lens, a transparent ball with varying radially symmetric refractive index \cite{E}.

\begin{figure}[h]
\begin{picture}(0,95)
\rput(5,0){
\rput(1,0.5){
\psline[linewidth=1.2pt](2,0)(0,0)(0,2)
\psline[linewidth=0.4pt,arrows=->,arrowscale=1.5,linecolor=red](2,0.5)(1,0)(0.2,0.4)
\psline[linewidth=0.4pt,arrows=->,arrowscale=1.5,linecolor=red](4,1.5)(2,0.5)
\psline[linewidth=0.4pt,arrows=->,arrowscale=1.5,linecolor=red](4,1.5)(1,0)(0,0.5)(3,2)
\psline[linewidth=0.4pt,arrows=->,arrowscale=1.5,linecolor=blue](2,1.9)(0,1.4)(2,0.9)
\pscircle[fillstyle=solid,linecolor=blue](2,2.2){0.25}
\rput(2,2.2){\blue 2}
\pscircle[fillstyle=solid,linecolor=red](4,1.1){0.25}
\rput(4,1.1){\red 1}
}
}
\end{picture}
\caption{Square corner: a \rr\ based on light reflection. Two incident light rays are shown: the ray 1 is retroreflected, while the ray 2 is not.}
\label{fig 2D practical retroreflectors}
\end{figure}
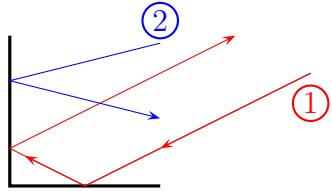

The most commonly used \rr\  based solely on {\it light reflection} is the so-called {cube corner} (its two-dimensional analogue, square corner, is shown in figure \ref{fig 2D practical retroreflectors}). Both cube and square corners are not perfect, however: a part of incoming light is reflected in a wrong direction. This is clearly seen in fig. \ref{fig 2D practical retroreflectors} for the square corner.

In what follows, only \rrs\ based on {\it light reflection} (or billiard \rrs) will be considered. To the best of our knowledge, no {\it perfect} billiard \rrs\ are known. However, as will be shown below, there exist \rrs\ which are {\it almost perfect}; more precisely, there exists a family of bodies $B_\ve, \ \ve > 0$ (which will be called an {\it asymptotically perfect \rr}) such that the portion of light reflected by $B_\ve$ in wrong directions goes to zero as $\ve \to 0$.

The main aim of this paper is twofold. First, bring together billiard type \rrs\ known by now. They form a small collection of four objects; the first, the second and the fourth one are \aprrs, and the third one is a \rr\ which is very close to perfect. The first three objects --- {\it mushroom}, {\it tube} and {\it \dopa} --- have already been published or submitted for publication \cite{Nonlinearity,tube,Gouveia Plakhov Torres AMC}. Note that the proof of retroreflectivity for the tube reduces to a quite nontrivial ergodic problem considered in \cite{tube}. The \dopa\ has been discovered and studied numerically \cite{Gouv disser,Gouveia Plakhov Torres AMC}. The fourth object --- {\it notched angle} --- is the new one. The second aim of the paper is to describe this shape and provide a proof of its retroreflectivity.

In section \ref{sect math definitions} we define basic mathematical notions that are used in the following sections \ref{sect collection} and \ref{sect notched angle}. The notions of perfect and \aprrs\ are introduced, and a quantity characterizing retroreflecting properties of a given body is determined. Also, in the two-dimensional case we introduce the notion of a hollow on the body boundary and describe billiard scattering in a hollow. In section \ref{sect collection} we present the collection of billiard \rrs\ and discuss and compare their properties. Finally, section \ref{sect notched angle} is devoted to the proof of retroreflectivity of notched angle, the fourth object in the collection.

\section{Mathematical preliminaries}
\label{sect math definitions}

Here we introduce basic notions and provide necessary information that will be used in the following sections.

Consider a connected set $B \subset \RRR^d$ with piecewise smooth boundary (in what follows such a set will be called {\it a body}), and consider the billiard in $\RRR^d \setminus B$. We shall denote by $x(t)$ the coordinate of a billiard particle at the moment $t$, by $v(t) = x'(t)$ its velocity, and by $v$ and $v^+$ the limits $v = \lim_{t\to -\infty} v(t), \ v^+ = \lim_{t\to +\infty} v(t)$, if they exist.

%\begin{opr}\label{o incident particle}
We say that a billiard particle is incident on $B$, if it moves freely prior to a moment $t_1$ and collides with $B$ at this moment. That is, the part of the trajectory $x(t), \ t < t_1$ is a half-line contained in $\RRR^d \setminus \bar B$ and $x(t_1) \in \pl B$.
%\end{opr}

\begin{opr}\label{o retroreflector n}
{\rm A body $B$ is called a {\it \prr}, if for almost all incident particles the asymptotic velocity at $t \to +\infty$ exists and is opposite to the asymptotic velocity at $t \to -\infty$; that is, $v^+ = -v$.}
\end{opr}

\begin{remark}\label{zam pathology}
{\rm Notice that the trajectory of some particles cannot be extended beyond a certain moment of time. This happens when the particle gets into a singular point of the boundary $\pl B$ or makes infinitely many reflections in a finite time. However, the set of such "pathological" particles has zero measure (see, e.g., \cite{Tab}) and will be excluded from our consideration.}
\end{remark}

\subsection{Unbounded bodies}
\label{subsec unbounded bodies}

The case of unbounded bodies is quite simple. Here we provide several examples of unbounded \prrs.
\vspace{2mm}

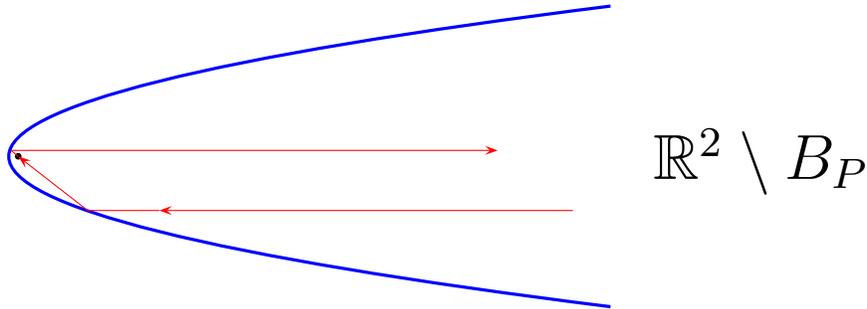
\begin{figure}[h]
\begin{picture}(0,130)
\rput(-0.6,-0.3){
\rput{90}(11,0.5){
\parabola[linewidth=1.2pt,linecolor=blue](0,0)(2,8)
\psdots[dotsize=2.5pt](2,7.875)%(2,7.875)(1.28,6.95)(2.08,7.97)
%\psline(1,6.575)(3,9.175)
 \psline[linewidth=0.4pt,linecolor=red,arrows=->,arrowscale=1.5](1.28,0.5)(1.28,6)
 \psline[linewidth=0.4pt,linecolor=red](1.28,6)(1.28,6.95)
 \psline[linewidth=0.4pt,linecolor=red,arrows=->,arrowscale=1.5](1.28,6.95)(2,7.875)
 \psline[linewidth=0.4pt,linecolor=red](2,7.875)(2.08,7.97)
 \psline[linewidth=0.4pt,linecolor=red,arrows=->,arrowscale=1.5](2.08,7.97)(2.08,1.5)
}
\rput(13,2.5){\scalebox{2}{$\RRR^2 \setminus B_P$}}

}
\end{picture}
\caption{Exterior of a parabola: an example of unbounded retroreflector with a unique velocity of incidence.}
\label{fig parabola}
\end{figure}

\hspace*{-6.5mm}{\bf Example 1.} $B = B_P$ is the exterior of a parabola in $\RRR^2$. There exists a unique velocity of incidence, which is parallel to the parabola axis. The initial and final velocities of any incident particle are mutually opposite, and the segment of the trajectory between the two consecutive reflections passes through the focus, as shown in figure \ref{fig parabola}.

\begin{remark}\label{zam retr parab}
{\rm If $B$ is the exterior of a parabola perturbed within a bounded set (that is, $B = B_P \vartriangle K$, with $K$ bounded), then $B$ is again a \prr. Indeed, any segment (or the extension of a segment) of a billiard trajectory within the parabola touches a confocal parabola with the same axis. The branches of this confocal parabola are co-directional or counter-directional with respect to the original parabola. This implies that the segments of an incident trajectory, when going away to the infinity, are becoming "straightened"{}, that is, more and more parallel to the parabola axis, and therefore $v^+ = -v$.}
\end{remark}

There also exist unbounded \rrs\ that admit a continuum of incidence velocities.
\vspace{2mm}

\hspace*{-6.5mm}{\bf Example 2.} Let $\RRR^d \setminus B$ be determined by the relations $x_1 > 0, \ldots, x_d > 0$ in an orthonormal reference system $x_1, \ldots, x_d$; then $B$ is a \prr.
\vspace{2mm}

Consider one more example.
\vspace{2mm}

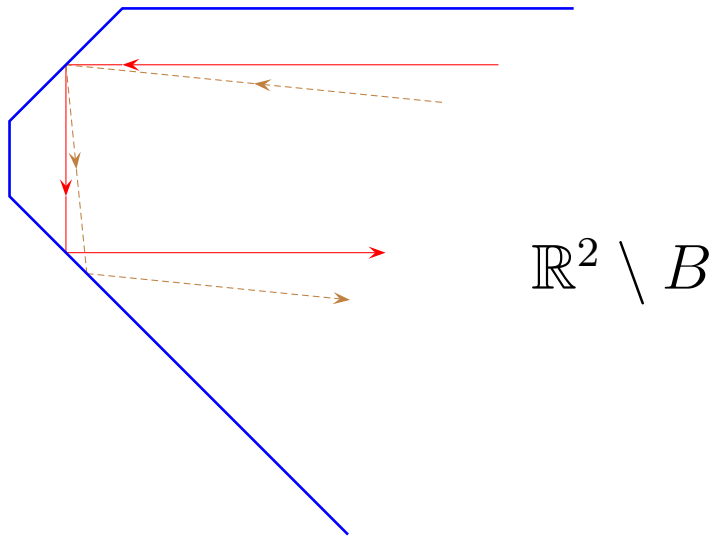
\begin{figure}[h]
\begin{picture}(0,200)
\rput(4,4.2){
\scalebox{0.5}{
\psline[linewidth=0.8pt,linecolor=red,arrows=->,arrowscale=3](13,3.5)(3,3.5)
\psline[linewidth=0.8pt,linecolor=red](3,3.5)(1.5,3.5)
\psline[linewidth=0.8pt,linecolor=red,arrows=->,arrowscale=3](1.5,3.5)(1.5,0)
\psline[linewidth=0.8pt,linecolor=red](1.5,0)(1.5,-1.5)
\psline[linewidth=0.8pt,linecolor=red,arrows=->,arrowscale=3](1.5,-1.5)(10,-1.5)
  \psline[linewidth=0.8pt,linecolor=brown,arrows=->,arrowscale=3,linestyle=dashed](11.5,2.5)(6.5,3)
  \psline[linewidth=0.8pt,linecolor=brown,linestyle=dashed](6.5,3)(1.5,3.5)
  \psline[linewidth=0.8pt,linecolor=brown,arrows=->,arrowscale=3,linestyle=dashed](1.5,3.5)(1.778,0.722)
  \psline[linewidth=0.8pt,linecolor=brown,linestyle=dashed](1.778,0.722)(2.056,-2.056)
  \psline[linewidth=0.8pt,linecolor=brown,arrows=->,arrowscale=3,linestyle=dashed](2.056,-2.056)(9.056,-2.756)
     \psline[linewidth=2pt,linecolor=blue](9,-9)(0,0)(0,2)(3,5)(15,5)
}
\rput(8,-1){\scalebox{2}{$\RRR^2 \setminus B$}}
}
\end{picture}
\caption{The two-dimensional unbounded \rr\ shown here is a convex polygon contained in an angle of size $\pi/4$, with all angles at its vertices being multiples of $\pi/4$. Two billiard trajectories in $\RRR^2 \setminus B$ are shown.}
\label{fig 2D unbounded retroreflectors}
\end{figure}

\hspace*{-6.5mm}{\bf Example 3.} Let the set $\mathbb C \setminus B$ in the complex plane $\mathbb{C} \sim \RRR^2$ be given by the relations Re$(e^{\frac{i\pi k}{2m}}z) > a_k$,\, $k = 0,\, 1, \ldots, 2m-1$, with $m \in \mathbb{N}$ and arbitrary constants $a_k$; then $B$ is a \prr; see figure \ref{fig 2D unbounded retroreflectors} for the case $m = 2$.

%\end{enumerate}

\subsection{Bounded bodies}
\label{subsec bounded bodies}

In what follows we restrict ourselves to the case of {\it bounded bodies}, which is more interesting both from mathematical viewpoint and for applications.

At present, no bounded \prrs\ are known. On the other hand, there exist families of bounded \rrs\ which are asymptotically perfect. The next two sections are devoted to description of and studying such families. Let us give exact definitions.

Consider a particle incident on $B$ that initially (prior to collisions with $B$) moves freely according to $x(t) = \xi + vt$, and denote by $v^+_B(\xi, v)$ its final velocity. The function $v_B^+$ is defined for all values $(\xi, v)$ such that the straight line $\xi + vt$,\, $t \in \RRR$ has nonzero intersection with $B$, except possibly for a set of zero measure.

Consider a convex body $C$ containing $B$ and define the measure $\mu_C$ on $\pl C \times S^{d-1}$ according to $d\mu_C(\xi, v) = \langle n(\xi),\, v \rangle_{\!-}\, d\xi\, dv$, where $n(\xi)$ is the outer normal to $\pl C$ at $\xi \in \pl C$,\, $d\xi$ and $dv$ are Lebesgue $(d-1)$-dimensional measures on $\pl C$ and $S^{d-1}$, respectively, $\langle \cdot\,, \cdot \rangle$ means scalar product, and $z_- = \max \{ 0,\, -z \}$ is the negative part of the real number $z$.

The mapping $T = T_{B,C} : (\xi, v) \mapsto (v, v_B^+(\xi, v))$ induces the push-forward measure $\nu_{B,C} = T_\# \mu_C$ on $(S^{d-1})^2$. One easily verifies (see \cite{UMN review article}) that $\nu_{B,C}$ does not depend on the ambient body $C$, and therefore one can just write $\nu_B$, omitting the subscript $C$. This measure admits a natural interpretation: it determines the (normalized) number of particles with initial and final velocities $v,\, v^+$ that have interacted with $B$ during a unit time interval.

\begin{opr}\label{o asympt perfect rr}
{\rm We say that $\nu$ is a {\it retroreflector measure}, if spt$\,\nu$ is contained in the subspace $\{ v^+ = -v \}$. A family of bounded bodies $B_\ve$,\, $\ve > 0$ is called an {\it \aprr}, if the measure $\nu_{B_\ve}$ weakly converges to a \rr\ measure as $\ve \to 0$.}
\end{opr}

\begin{remark}\label{zam bounded rr}
{\rm From definition \ref{o retroreflector n} it follows that a bounded body is a \prr\ {\it iff} $\nu_B$ is a \rr\ measure.}
\end{remark}

In the two-dimensional case one easily calculates the full measure $\nu_B((S^1)^2)$. Take $C =$ Conv$\,B$; then, introducing the natural parameter $\xi \in [0,\, |\pl C|]$ on $\pl C$ and denoting by $\vphi \in \interval$ the angle (counted counterclockwise) from $-n(\xi)$ to $v$, one gets
$$
\nu_B((S^1)^2) = \mu_C(\pl C \times S^1) = \int_0^{|\pl C|} d\xi \int_{-\pi/2}^{\pi/2} \cos\vphi\, d\vphi = 2|\pl C|.
$$

\subsection{Resistance}\label{subsec resistance}

Here we introduce a functional on the set of bounded bodies that indicates how close the billiard scattering by the body is to the retroreflector scattering. This functional is called {\it normalized resistance}, a quantity that has mechanical interpretation going back to Newton's problem of minimal resistance \cite{N}. We believe that it serves as a natural measure of "retroreflectivity".

The force of resistance of the body $B$ to a parallel flow of particles at the velocity $v$ equals
 \begin{equation}\label{resist B,v}
R(B,v) = \int_{v^\perp} (v - v^+_B(\xi, v))\, d\xi,
 \end{equation}
where $v^\perp$ is the orthogonal complement to the one-dimensional subspace $\{ v \}$. (We suppose that the flow has unit density.) The expression (\ref{resist B,v}) is defined for almost all $v \in S^{d-1}$. The component of the resistance force along the flow direction equals $\langle R(B,v),\, v \rangle$.

Suppose that the velocity of the flow $v$ is taken at random and uniformly in $S^{d-1}$; then the mathematical expectation of the resistance along the flow equals $\mathbb{E}\langle R(B,v),\, v \rangle = c\, R(B)$, where $c = 1/|S^{d-1}|$ and
 \begin{equation}\label{mean resist}
R(B) = \int_{S^{d-1}} \langle R(B,v),\, v \rangle\, dv.
 \end{equation}

Let $C$ be a convex body containing $B$. Taking into account the invariance of $v^+_B$ relative to translations along $v$,\, $ v^+_B(\xi, v) =  v^+_B(\xi + vt, v)$ and making a change of variables, the integral $R(B)$ can be transformed to the form
 %\begin{equation}\label{resist}
 $$
R(B) =  \int_{\pl C \times S^{d-1}} \langle v - v^+_B(\xi, v),\, v \rangle \langle v,\, n(\xi) \rangle_{\!-}\, d\xi\, dv  \quad \quad \quad \quad \quad \quad \quad \quad \quad \quad \quad
 $$
  %\end{equation}
 \begin{equation}\label{resist}
 \quad \quad \quad \quad \quad \quad \quad \quad \quad \quad \quad \quad \quad \quad \quad = \int_{\pl C \times S^{d-1}} \langle v - v^+_B(\xi, v),\, v \rangle\, d\mu_C(\xi, v).
 \end{equation}
  %The integral in (\ref{resist}) does not depend on the choice of ambient convex body $C$. In particular, one can take $C = $ Conv$B$.
Using the definition of $\nu_B$ and making one more change of variables, one gets
 \begin{equation}\label{resist nu}
R(B) = \int_{(S^{d-1})^2} \left( 1 - \langle v,\, v^+ \rangle \right) d\nu_B(v,v^+).
 \end{equation}

\begin{remark}\label{zam R mech interpretation}
{\rm Let us mention another mechanical interpretation of the quantity $R(B)$. Suppose that the body $B$ translates through a medium of resting particles and at the same time slowly and chaotically rotates (somersaults), so that in a reference system connected with the body the vector of translational velocity runs $S^{d-1}$ chaotically and uniformly. Then the mean value of resistance during a long period of time approaches $R(B)$ when the length of the period goes to infinity.}
\end{remark}

Let us additionally define the mean resistance of the body under the so-called {\it diffuse scattering}, where each incident particle, after hitting the body, completely loses its initial velocity and remains near $\pl B$ forever. The formula for the diffuse resistance, $D(B)$, is similar to the above formula (\ref{resist nu}) for the elastic resistance. The difference is that the normalized momentum transmitted by a particle to the body is always equal to 1, and therefore, the integrand $1 - \langle v,\, v^+ \rangle$ in (\ref{resist nu}) should be substituted with 1. The resulting formula is
 \begin{equation}\label{diffuse resist}
D(B) = \int_{(S^{d-1})^2} d\nu_B(v,v^+) = \nu_B((S^{d-1})^2).
 \end{equation}

Notice that the following inequality always holds
$$
R(B) \le 2D(B);
$$
besides, if $B$ is a hypothetical \rr, this inequality turns into the equality $R(B) = 2D(B)$.

\begin{remark}\label{zam diffuse scattering}
{\rm The notion of diffuse scattering has a strong physical motivation originating, in particular, from space aerodynamics. The interaction of artificial satellites on low Earth orbits with the rarefied atmosphere is considered to be mainly diffuse by some researches (see, e.g., \cite{MM}). Some others (\cite{W,IY,BSE}) prefer to use Maxwellian representation of interaction as a linear combination of elastic scattering and diffuse one. In the latter case the resistance equals $\al D(B) + (1 - \al) R(B)$, where $\al$ is the so-called accommodation coefficient.}
\end{remark}

Let us calculate $R(B)$ and $D(B)$ in the case where $B$ is convex. Using (\ref{resist}) and taking into account the formula of elastic scattering $v^+ = v - 2\langle v,\, n \rangle n$, one gets
$$
\frac{1}{|\pl B|}\, R(B) = \frac{1}{|\pl B|} \int_{\pl B \times S^{d-1}} 2\langle v,\, n(\xi) \rangle^2 \langle v,\, n(\xi) \rangle_-\, d\xi\, dv = \int_{S^{d-1}} 2\langle v,\, n \rangle_-^3\, dv =
$$
$$
= |S^{d-2}| \int_0^{\pi/2} 2\cos^3 \vphi\, \sin^{d-2} \vphi\, d\vphi = \frac{4}{d+1}\, \frac{\pi^{\frac{d-1}{2}}}{\Gam(\frac{d+1}{2})},
$$
where $n$ is an arbitrary unit vector, and similarly,
$$
\frac{1}{|\pl B|}\, D(B) = \int_{S^{d-1}} \langle v,\, n \rangle_{\!-}\, dv = |S^{d-2}| \int_0^{\pi/2} \cos\vphi\, \sin^{d-2} \vphi\, d\vphi = \frac{\pi^{\frac{d-1}{2}}}{\Gam(\frac{d+1}{2})}.
$$
Therefore one has
$$
\frac{R(B)}{D(B)} = \frac{4}{d+1}.
$$
In particular, in the three-dimensional case one gets the equality $R(B) = D(B)$; that is, the elastic resistance of convex bodies is equal to the diffuse one.

Define the {\it normalized mean resistance} of the body as follows:
 \begin{equation}\label{specific resist}
r(B) = \frac{R(B)}{2D(B)}.
 \end{equation}
It has the following useful properties.

\begin{enumerate}

\item $0 \le r(B) \le 1$.

\item If $B$ is convex then $r(B) = 2/(d+1)$; in particular, $r(B) = 2/3$ for $d=2$ and $r(B) = 1/2$ for $d=3$.

\item $\sup_B r(B) = 1$ in any dimension.

\item The infimum of $r$ depends on the dimension $d$.\\
In the case $d = 2$,\, $\inf_B r(B) = 0.6585...$ (see \cite{ARMA}).\\
In the case $d \ge 3$ only estimates are known. In particular, if $d = 3$ then $\inf_B r(B) < 0.4848$ (see \cite{UMN review article}).

\item  If $B_\ve$ is an \aprr\ then $\lim_{\ve\to 0} r(B_\ve) = 1$.

\end{enumerate}

The property 3 is a consequence of existence, in any dimension, of \aprrs\ (see subsection \ref{subsec mushroom}).

\begin{remark}\label{zam norm res}
{\rm The value $r(B)$ is proportional to the (elastic) resistance of $B$ divided by the number of particles that have interacted with $B$ during a unit time interval. It can also be interpreted as the mathematical expectation of the longitudinal component of the momentum transmitted to the body by a randomly chosen incident particle of mass $1/2$, that is, $r(B) = \frac 12\, \mathbb{E} \langle v-v^+,\, v \rangle$}.
\end{remark}

\subsection{Hollow}\label{subsec hollow and scattering}

\begin{figure}[h]
\begin{picture}(0,230)
\rput(4,0){
\scalebox{0.8}{
\psecurve[linewidth=1pt,fillcolor=yellow,fillstyle=solid](1.5,4.75)(0.25,6.5)(0.8,8.5)(2.5,9.75)(4.5,9.5)(4.8,8)(5.5,7.1)(7.5,7)(9,7.3)(10.3,5.5)(8.7,4.3)
(7.7,3)(7.5,1)(6,0.25)(4.5,0.6)(3.25,2)(3.5,3.25)(3,4.3)(1.3,4.5)(0.25,6.5)(0.8,8.5)
\psecurve[linewidth=1pt,fillcolor=white,fillstyle=solid](4,4.6)(5,4.15)(6,4.3)(6,5.1)(5,5.4)(4.3,6.2)(3.5,7.3)(2.1,7)(2.3,6)(3,5.3)(4,4.6)(5,4.15)(6,4.3)
\psline[linestyle=dashed,linewidth=0.5pt](0.8,5)(3.55,1.4)
\psline[linestyle=dashed,linewidth=0.5pt](4.22,9.71)(9.15,7.25)
\psline[linestyle=dashed,linewidth=0.5pt](7.52,1.02)(10.2,5.23)
\rput(6,2.8){\scalebox{4}{$B$}}
\rput(2.8,3.6){\scalebox{1.5}{$\Om_1$}}
\rput(6,7.8){\scalebox{1.5}{$\Om_2$}}
\rput(8.4,3.3){\scalebox{1.5}{$\Om_3$}}
\rput(3.6,5.9){\scalebox{1.5}{$\Om_4$}}
\rput(1.8,3.1){\scalebox{1.5}{$I_1$}}
\rput(7,8.8){\scalebox{1.5}{$I_2$}}
\rput(9.25,3){\scalebox{1.5}{$I_3$}}
\rput(10.5,6.5){\scalebox{1.5}{$I_0$}}
\rput(0.45,8.75){\scalebox{1.5}{$I_0$}}
\rput(4.5,0.1){\scalebox{1.5}{$I_0$}}
}
}
\end{picture}
\caption{A body $B$ and the corresponding hollows.}
\label{fig set with hollows}
\end{figure}
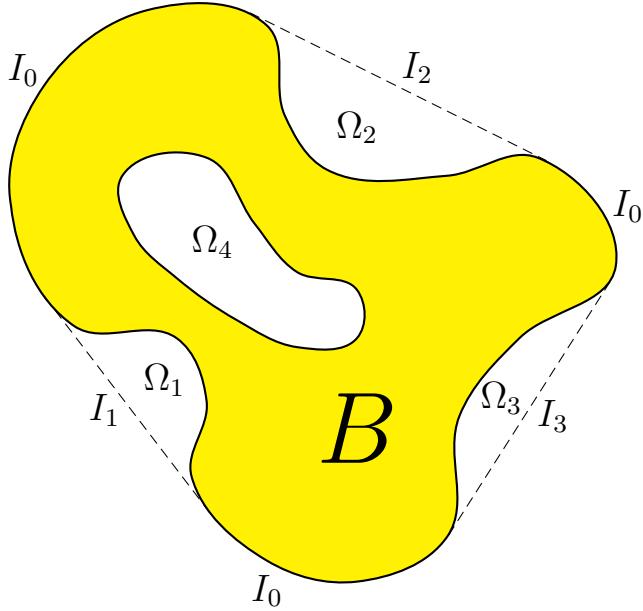

Here we consider the two-dimensional case, $d = 2$. Take a bounded body $B$ and represent each of the sets $\text{Conv}B \setminus B$ and $\pl(\text{Conv}B) \setminus \pl B$ as the union of its connected components,
$$
\text{Conv}B \setminus B = \cup_{i\ne 0} \Om_i, \quad \pl(\text{Conv}B) \setminus \pl B = \cup_{i\ne 0} I_i;
$$
in both cases the set of indices $i$ is finite or countable, and each set $I_i$ is an interval contained in $\Om_i$; see fig. \ref{fig set with hollows}. (Notice that $B$ is not necessarily simply connected, and so, there may exist sets $\Om_i$ that are entirely contained in $B$ and therefore do not contain any interval $I_j$.) Denote by $I_0$ the convex part of the boundary $\pl B$,\, $I_0 = \pl(\text{Conv}B) \cap \pl B$; thus, one has
$$
\pl(\text{Conv}B) = \cup_{i} I_i.
$$

\begin{opr}\label{o hollow}
{\rm Any pair of sets $(\Om_i, I_i)$ that appears in the above construction applied to a bounded body $B$ is called a {\it hollow}. The interval $I_i$ is called the} {opening of the hollow}.
\end{opr}

Any hollow $(\Om,I)$ has the following properties.

 \begin{quote}
(i) $\Om$ is a bounded simply connected set with piecewise smooth boundary.

(ii) $I$ is an interval contained in $\pl\Om$.

(iii) $\Om$ is situated on one side of the straight line containing $I$.

(iv) The intersection of $\Om$ with this line coincides with $I$.
 \end{quote}

Inversely, any pair $(\Om,I)$ satisfying the conditions (i)--(iv) is a hollow.

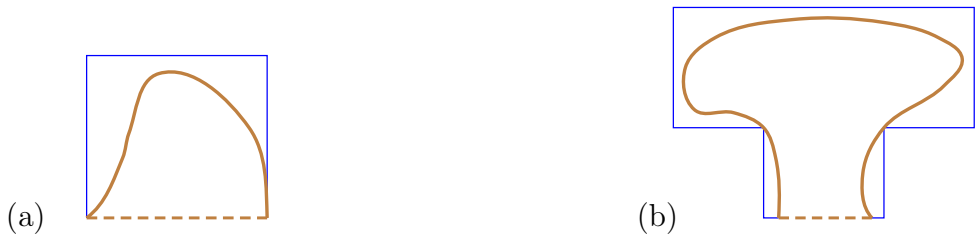
\begin{figure}[h]
\begin{picture}(0,100)
\rput(2,0.3){
\rput(-0.7,0){(a)}
\scalebox{0.8}{
\psline[linewidth=0.6pt,linecolor=blue](0,0)(0,2.7)(3,2.7)(3,0)
\psline[linestyle=dashed,linewidth=1.4pt,linecolor=brown](0,0)(3,0)
\psecurve[linewidth=1.6pt,linecolor=brown](-1,-0.25)(0,0)(0.6,1)(0.7,1.4)(1.2,2.4)(2.7,1.5)(3,0)(2.8,-1.5)
}
}
\rput(11,0.3){
\rput(-1.3,0){(b)}
\scalebox{0.8}{
\psline[linewidth=0.6pt,linecolor=blue](0.25,0)(0,0)(0,1.5)(-1.5,1.5)(-1.5,3.5)(3.5,3.5)(3.5,1.5)(2,1.5)(2,0)(1.8,0)
%\psline[linestyle=dashed,linewidth=0.6pt,linecolor=brown](0,0)(2,0)
\psline[linestyle=dashed,linewidth=1.4pt,linecolor=brown](0.25,0)(1.8,0)
\psecurve[linewidth=1.6pt,linecolor=brown](0.1,-1)(0.25,0)(0.2,1)(0,1.5)(-0.5,1.75)(-1.15,1.8)(-1.3,2.5)(-0.8,3)(0.5,3.3)
(2,3.25)(2.9,3)(3.3,2.65)(3,2.2)(2,1.5)(1.65,0.75)(1.8,0)(2.5,-0.2)
}
}
\end{picture}
\caption{(a) A convenient hollow. (b) An inconvenient hollow.}
\label{fig convenient hollows}
\end{figure}

\begin{opr}\label{o convenient hollow}
{\rm A hollow $(\Om, I)$ is called {\it convenient}, if the orthogonal projection of $\Om$ on the line containing $I$ coincides with $I$. Otherwise, it is called {\it inconvenient}. See figures \ref{fig convenient hollows}a and \ref{fig convenient hollows}b for examples of convenient and inconvenient hollows}.
\end{opr}

An incident particle may hit the body in the convex part of its boundary and then go away. Otherwise, it gets into a hollow through its opening, makes there several reflections, and then escapes the hollow through the opening and goes away. It is helpful to define the {\it measures generated by hollows} and the measure generated by the convex part of the boundary, and then represent $\nu_B$ as a weighted sum of these measures.

\begin{figure}[h]
\begin{picture}(0,110)
\rput(0.3,-2.5){
\scalebox{1.13}{
\rput{31}(3.2,-4){
\scalebox{0.7}{
 \psecurve[linewidth=1.1pt,linecolor=brown]%(3,8)(2.5,6.5)(1.7,5.4)(0.5,4.7)(0.1,4)(0,2.5)(1,1)(3,0.2)(6,0.2)(9,1.25)
 (11.4,3)(12.3,4.5)(12.4,5.35)(12,5.5)(11,5.2)(9.5,5.2)(8.3,6)(7.9,7)(8,8)(7.8,8)(7.6,8.1)
 \psline[linestyle=dashed,linewidth=0.6pt,linecolor=green](12.3,5.45)(8,8)
    \psarc[linewidth=1pt,linecolor=brown](8.3,7.82){1.1}{59.331}{90}
    \psarc[linewidth=1pt,linecolor=brown](11.76,5.77){1}{36.87}{59.331}
    \psarc[linewidth=1pt,linecolor=brown](11.76,5.77){1.1}{36.87}{59.331}
 \psline[linewidth=1.2pt,arrows=->,arrowscale=1.5,linecolor=red](8.3,9.5)(8.3,6)(11,5.2)(13.4,7)%(12.6,6.4)
 \psline[linewidth=1.2pt,arrows=->,arrowscale=1.5,linecolor=red](8.3,9.7)(8.3,8.5)
 \psline[linewidth=1pt,arrows=->,arrowscale=1.5,linecolor=blue](11.76,5.77)(12.78,7.49)
 \psline[linewidth=1pt,arrows=->,arrowscale=1.5,linecolor=blue](8.3,7.82)(9.32,9.54)
 \psdots[dotsize=3pt](8.3,7.82)(11.76,5.77)
    }
    }
\scalebox{0.7}{
     \rput(6.5,6.3){\Large $v$}
  \rput(7.23,6.6){\large $\vphi$}
   \rput(11.9,6.75){$\vphi^+$}
 \rput(12.8,6.9){\Large $v^+$}
  \rput(11.7,7.56){\Large $n_i$}
    \rput(7.8,7.56){\Large $n_i$}
 \rput(7.8,5.67){\Large $\xi$}
 \rput(11.25,5.5){\Large $\xi^+$}
        \rput(9.9,5.7){\Large $I_i$}
        \rput(9.1,4.6){\Large $\Om_i$}
        }
        }
        }
    \end{picture}
\caption{Billiard scattering in a hollow. Here one has $\vphi > 0$ and $\vphi^+ < 0$.}
\label{fig scat on hollow}
\end{figure}
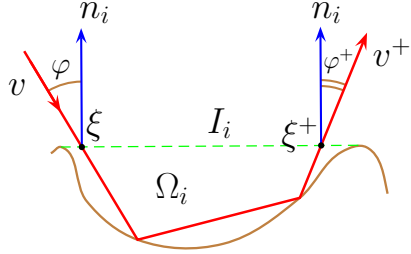

Consider a hollow $(\Om_i, I_i)$, denote by $n_i$ the outer normal to $\Om_i$ at an arbitrary point of $I_i$, and introduce a uniform coordinate $\xi \in [0,\, 1]$ on $I_i$ varying from 0 at one endpoint of $I_i$ to 1 at the other one. For a particle that gets into the hollow at the velocity $v$, makes there several (maybe one) reflections, and then gets out at the velocity $v^+$, fix the coordinate $\xi$ of the first intersection with $I_i$ (when the particle gets "in"), denote by $\vphi$ the angle between the vectors $n_i$ and $-v$, and denote by $\vphi^+ = \vphi^+_{\Om_i,I_i}(\xi,\vphi)$ the angle between the vectors $n_i$ and $v^+$; see fig. \ref{fig scat on hollow}.
The angles are counted counterclockwise from $\pm n_i$ to $v$ or $v^+$; both angles belong to $\interval\!\! \mod\! 2\pi$. Define the probability measure $\mu$ on $[0,\, 1] \times \interval$ by
$$
d\mu(\xi,\vphi) = \frac 12 \cos\vphi\, d\xi\, d\vphi,
$$
where both $d\xi$ and $d\vphi$ denote one-dimensional Lebesgue measure.

The mapping $\TTT_i : (\xi, \vphi) \mapsto (\vphi, \vphi^+_{\Om_i,I_i}(\xi, \vphi))$ induces the push-forward measure  $\eta_{\Om_i,I_i} := \TTT_i^\# \mu$ on the square $\Box := \interval \times \interval$. Thus, one has
$$
\eta_{\Om_i,I_i}(A) = \mu\left(\left\{ (\xi,\vphi) : (\vphi, \vphi^+_{\Om_i,I_i}(\xi,\vphi)) \in A \right\}\right)
$$
for any Borel set $A \subset \Box$.

The probability measure $\eta_{\Om_i,I_i}$ is called the {\it measure generated by the hollow} $(\Om_i,I_i)$.

Notice that geometrically similar hollows generate identical measures.

Next, define the measure $\eta_{I_0}$ on $\Box \times I_0$ with the density $\frac{1}{2|I_0|}\, \cos\vphi\, \del(\vphi + \vphi^+)$, which will be called the {\it measure generated by the convex part of the boundary}, and the measure $\eta_0$ on $\Box$ with the density $\frac 12\, \cos\vphi\, \del(\vphi + \vphi^+)$.

Let $v_n(\vphi)$ be the vector obtained by rotating the vector $n$ counterclockwise by the angle $\vphi$, and let $n(\xi)$ be the outer normal to $B$ at $\xi \in \pl B$. The mapping $\s_i : (\vphi,\vphi^+) \mapsto (v_{-n_i}(\vphi),\, v_{n_i}(\vphi^+))$ induces the push-forward probability measure $\nu_{\Om_i,I_i} = \s_i^\# \eta_{\Om_i,I_i}$ on $(S^1)^2 = \mathbb{T}^2$, and the mapping $\s_0 : (\vphi,\vphi^+,\xi) \mapsto (v_{-n(\xi)}(\vphi),\, v_{n(\xi)}(\vphi^+))$ induces the push-forward probability measure $\nu_{I_0} = \s_0^\# \eta_{I_0}$ on $\mathbb T^2$. The measures $\nu_{\Om_i,I_i}$ and $\nu_{I_0}$ will also be called the measures generated by the hollows and the measure generated by the convex part of the boundary, respectively.

\begin{remark}\label{zam rr measure}
{\rm Consider the probability measure $\eta_\star$ on $\Box$ with the density $\frac 12\, \cos\vphi\, \del(\vphi-\vphi^+)$. Its push-forward measure $\s_i^\# \eta_\star$, for any $i$, is a \rr\ measure on $\mathbb{T}^2$. For this reason, $\eta_\star$ will also be called a \rr\ measure}.
\end{remark}

\begin{opr}\label{o as rr hollows}
{\rm A family of hollows $(\Om_\ve, I_\ve)$ is called {\it asymptotically retroreflecting}, if $\eta_{\Om_\ve,I_\ve}$ weakly converges to $\eta_\star$.}
\end{opr}

The measure $\nu_B$ can be represented as
$$
\nu_B = |I_0|\, \nu_{I_0} + \sum_{i\ne 0} |I_i|\, \nu_{\Om_i,I_i},
$$
and the functionals $R(B)$ and $D(B)$ then take the form
\begin{equation}\label{R(B)}
R(B) = |I_0| \int\!\!\!\!\int_{\mathbb{T}^2} (1 - \langle v,\, v^+ \rangle)\, d\nu_{I_0}(v,v^+) + \sum_{i\ne 0} |I_i| \int\!\!\!\!\int_{\mathbb{T}^2} (1 - \langle v,\, v^+ \rangle)\, d\nu_{\Om_i,I_i}(v,v^+),
\end{equation}
\begin{equation}\label{D(B)}
D(B) = \nu_B(\mathbb{T}^2) = |I_0| + \sum_{i\ne 0} |I_i| = |\pl(\text{Conv}\,B)|.
\end{equation}
Using (\ref{R(B)}) and the relation between the measures $\eta_{\Om_i,I_i}$,\, $\eta_{0}$ and the measures $\nu_{\Om_i,I_i}$,\, $\nu_{0}$, and taking into account that $\langle v,\, v^+ \rangle = -\cos(\vphi - \vphi^+)$, one gets
$$
R(B) = |I_0| \int\!\!\!\int_{\Box} (1 + \cos(\vphi-\vphi^+))\, d\eta_{0}(\vphi,\vphi^+) + \hspace*{65mm}
$$
\begin{equation}\label{R(B)eta}
\hspace*{50mm} + \sum_{i\ne 0} |I_i| \int\!\!\!\int_{\Box} (1 + \cos(\vphi-\vphi^+))\, d\eta_{\Om_i,I_i}(\vphi,\vphi^+).
\end{equation}
Denote $c_i = |I_i|/|\pl(\text{Conv}B)|$,\, $\sum c_i = 1$ and define the functional
$$
\FFF(\eta) = \frac 12 \int\!\!\!\int_\Box (1 + \cos(\vphi-\vphi^+))\, d\eta(\vphi,\vphi^+).
$$
One easily calculates that $\FFF(\eta_0) = 2/3$ and $\FFF(\eta_\star) = 1$. Then, using (\ref{specific resist}), (\ref{R(B)eta}) and (\ref{D(B)}), one obtains
\begin{equation}\label{r(B)}
r(B) = \frac 23\, c_0 + \sum_{i\ne 0} c_i \FFF(\eta_{\Om_i,I_i}).
\end{equation}

The formula (\ref{r(B)}) suggests a strategy of constructing \aprrs. First, find an asymptotically retroreflecting family of hollows $(\Om_\ve, I_\ve)$; that is, $\lim_{\ve\to 0} \FFF(\eta_{\Om_\ve,I_\ve}) = 1$. Then find a family of bodies $B_\ve$ with all hollows on their boundary similar to $(\Om_\ve, I_\ve)$ and such that the relative length of the convex part of $\pl B_\ve$ goes to zero, $\lim_{\ve\to 0} c_0^\ve = 0$, and the sequence of convex hulls Conv$\,B_\ve$ converges to a fixed convex body as $\ve\to 0$. In this case one has
$$
\lim_{\ve\to 0} r(B_\ve) = \lim_{\ve\to 0} \left( \frac 23 c_0^\ve + (1 - c_0^\ve) \FFF(\eta_{\Om_\ve,I_\ve})\right) = 1,
$$
and therefore, the family $B_\ve$ is an \aprr.

If all the hollows are convenient (see fig. \ref{fig convenient hollows}a), then one can find bodies $B_\ve$ with identical hollows. If the hollows are not convenient (see fig. \ref{fig convenient hollows}b), then each body $B_\ve$ must contain, on its boundary, a hierarchy of hollows of different sizes.

\subsection{Semi-retroreflecting hollows}

\begin{figure}[h]
\begin{picture}(0,95)
\rput(2,0){
\psline[linecolor=brown](0,0)(4,0)(4,1)(0,1)
\psline[linewidth=0.3pt]{<->}(04.3,0)(4.3,1)
\rput(4.5,0.5){$\ve$}
\psline[linewidth=0.3pt]{<-}(0,1.4)(1.8,1.4)
\psline[linewidth=0.3pt]{->}(2.2,1.4)(4,1.4)
\rput(2,1.4){$1$}
\psline[linestyle=dotted,linewidth=1pt,linecolor=brown](0,0)(0,1)
\rput(-1.2,0.3){(a)}
}
\rput(13,-0.2){
\psline[linecolor=brown](-0.7,0)(0,3.5)(0.7,0)
\psarc[linecolor=blue](0,3.5){0.3}{259}{281}
\psarc[linecolor=blue](0,3.5){0.4}{259}{281}
\rput(-0.25,3.2){$\ve$}
\psline[linestyle=dotted,linewidth=1pt,linecolor=brown](-0.7,0)(0.7,0)
\rput(-1.6,0.3){(b)}
}
\end{picture}
\caption{A rectangular hollow (a) and a triangular hollow (b).}
\label{fig simple shapes}
\end{figure}
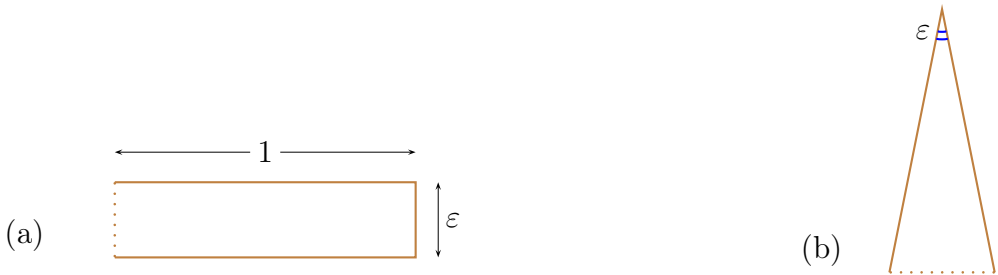

Let us mention two special kinds of hollows, a rectangle and a triangle, as shown in figure \ref{fig simple shapes}. The ratio of the width to the height of the rectangle equals $\ve$. The triangle is isosceles, and the angle at the apex equals $\ve$.
Denote by $\nu_\sqcup^\ve$ and $\nu_\vee^\ve$ the measures generated by the rectangle and the triangle, respectively.

\begin{utv}\label{utv semi}
Both $\nu_\sqcup^\ve$ and $\nu_\vee^\ve$ weakly converge to $\frac{1}{2}\, (\eta_0 + \eta_\star)$ as $\ve \to 0$.
\end{utv}

The proof of this statement is not difficult, but a little bit lengthy, and therefore is put in the appendix. The statement implies that both functionals, $\FFF(\nu^\ve_{\sqcup})$ and $\FFF(\nu^\ve_{\vee})$, converge to $5/6$. Note also that the measures $\nu_\sqcup^\ve$ and $\nu_\vee^\ve$ do not converge in norm.

Both the shapes are, so to say, semi-retroreflecting: nearly one half of the particles is reflected according to the elastic law $\vphi^+ = -\vphi$, and the other half, according to the \rr\ one $\vphi^+ = \vphi$. However, these shapes served as starting points for developing {\it true} \rrs: rectangular tube (subsection \ref{subsec tube}) and notched angle (subsection \ref{subsec notch ang} and section \ref{sect notched angle}).

\section{Collection of retroreflectors}\label{sect collection}

For each of the \aprrs\ proposed below, we first define the generating hollow $(\Om_\ve, I_\ve)$, and then construct the body $B_\ve$ formed by copies of this hollow.

\subsection{Mushroom}\label{subsec mushroom}

\begin{figure}[!hbp]
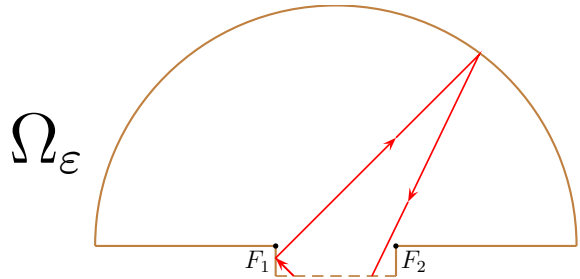

\pspicture*(-0.75cm,-0.5cm)(12cm,3.2cm)
 \rput(7.1,0){
 \scalebox{0.8}{
  \psarc[linewidth=1pt,linecolor=brown](0,0){4}{0}{180}
  \psline[linewidth=1pt,linecolor=brown](-4,0)(-1,0)(-1,-0.5)
  \psline[linewidth=1pt,linecolor=brown](4,0)(1,0)(1,-0.5)
  \psline[linewidth=0.8pt,linestyle=dashed,linecolor=brown](-1,-0.5)(1,-0.5)
  \psline[arrows=->,arrowscale=1.5,linecolor=red](-0.7,-0.5)(-1,-0.2)(1,1.8)
      \psline[arrows=->,arrowscale=1.5,linecolor=red](-0.7,-0.5)(-1,-0.2)
  \psline[arrows=->,arrowscale=1.5,linecolor=red](1,1.8)(2.38,3.18)(1.2,0.7333)
  \psline[linecolor=red](1.2,0.7333)(0.6,-0.5)
  \psdots[dotsize=2.5pt](-1,0)(1,0)
  \rput(-1.3,-0.25){$F_1$}
  \rput(1.3,-0.25){$F_2$}
  \rput(-4.8,1.7){\scalebox{2.8}{$\Om_\ve$}}
  }}
\endpspicture
\caption{Mushroom.}
\label{fig mushroom}
\end{figure}

The {\it mushroom} is the union of the upper semi-ellipse $\frac{x_1^2}{1+\ve^2} + x_2^2 = 1$,\, $x_2 \ge 0$  and the rectangle $-\ve \le x_1 \le \ve$,\, $-\ve^2 \le x_2 \le 0$ (see fig. \ref{fig mushroom}).
Its opening is the base of the mushroom stem, that is, the interval $[-\ve,\, \ve] \times \{-\ve^2\}$. The foci $F_1$ and $F_2$ of the ellipse are vertices of the rectangle, the width of the rectangle equals the focal distance $|F_1 F_2| = 2\ve$, and the heights, $\ve^2$.

Recall a remarkable property of the billiard in an ellipse. Any particle emanated from a focus, makes a reflection from the ellipse and then gets into the other focus. This implies that any particle that intersects the segment $F_1F_2$ in the direction "up"{}, after a reflection from the upper semi-ellipse will intersect this segment again, this time in the direction "down". Therefore, all particles getting into the mushroom through the opening, except for a portion $O(\ve)$, will make exactly one reflection and then get out, without hitting the mushroom stem. (The billiard trajectory depicted in figure \ref{fig mushroom} hits the stem, and therefore is exceptional.) For the non-exceptional particles the difference between the initial and final angle equals $\vphi - \vphi^+ = O(\ve)$. This simple observation leads to the following theorem.

\begin{theorem}\label{t mushroom}
The measure generated by mushroom weakly converges to $\eta_\star$ as $\ve \to 0$.
\end{theorem}

This theorem means that the mushroom is an asymptotically retroreflecting hollow. The mushroom and mushroom "seedlings"{} are discussed in \cite{Nonlinearity} in more detail.

\begin{remark}\label{bunimovich mushroom}
{\rm Notice that the mushroom was first introduced in billiard theory by Bunimovich as an example of dynamical system with divided phase space \cite{Bun}.}
\end{remark}

Let us describe some properties of the mushroom.
\vspace{1mm}

1. The mushroom is an {\it inconvenient} hollow. Therefore the resulting body (\aprr) contains a hierarchy of mushrooms of different sizes; see fig. \ref{fig Four bodies}(a).
\vspace{1mm}

%\begin{figure} \begin{center} \includegraphics[width=0.48\columnwidth]{cogumelosCirc} \end{center} \caption{A body with mushroom-shaped hollows.} \label{fig mushroom body} \end{figure}

2. The difference $\vphi - \vphi^+$ is always nonzero; this means that the mushroom measure converges to $\eta_\star$ {\it weakly}, but {\it not in norm}.
\vspace{1mm}

3. If the semi-ellipse is substituted with a semicircle then the resulting hollow (which is also called mushroom) will also be asymptotically retroreflecting. This modified construction can be generalized to any dimension; that is, there exist multidimensional \aprrs\ with mushroom-shaped hollows (for a more detailed description, see \cite{inverting UMN-06}).

4. Most incident particles make {\it exactly one reflection}. This means that the portion of incident particles making one reflection tends to 1 as $\ve \to 0$.

\subsection{Tube}\label{subsec tube}

The {\it tube} is a rectangle of width $a$ and height 1 with two rows of rectangles of smaller size $\del \times \ve$ taken away (see fig. \ref{fig tube}).
\begin{figure}[h]
\begin{picture}(0,90)
\rput(2.5,0.3){
\psline[linecolor=brown](0,0)(1.95,0)(1.95,0.3)(2.05,0.3)(2.05,0)(3.95,0)(3.95,0.3)(4.05,0.3)(4.05,0)
(5.95,0)(5.95,0.3)(6.05,0.3)(6.05,0)(7.95,0)(7.95,0.3)(8.05,0.3)(8.05,0)(9,0)(9,2)
\psline[linecolor=brown](0,2)(1.95,2)(1.95,1.7)(2.05,1.7)(2.05,2)(3.95,2)(3.95,1.7)(4.05,1.7)(4.05,2)
(5.95,2)(5.95,1.7)(6.05,1.73)(6.05,2)(7.95,2)(7.95,1.7)(8.05,1.7)(8.05,2)(9,2)
\psline[linewidth=0.3pt,linecolor=blue]{<-}(9.4,0)(9.4,0.75)
  \psline[linewidth=0.3pt,linecolor=blue]{<-}(0,2.4)(4.2,2.4)
  \rput(4.5,2.4){$a$}
  \psline[linewidth=0.3pt,linecolor=blue]{->}(4.8,2.4)(9,2.4)
  %%%%%%%%%%%%%%%%%%%%%%%%%%%%%%%%%%%%%%%%%%%
  \psline[linewidth=0.3pt,linecolor=blue]{->}(5.65,0.5)(5.95,0.5)
  \rput(6,0.8){$\del$}
  \psline[linewidth=0.3pt,linecolor=blue]{<-}(6.05,0.5)(6.35,0.5)
\psline[linewidth=0.3pt,linecolor=blue]{->}(9.4,1.25)(9.4,2)
\rput(9.4,1){1}
\psline[linewidth=0.3pt,linecolor=blue]{<-}(2.05,-0.3)(2.75,-0.3)
\psline[linewidth=0.3pt,linecolor=blue]{->}(3.25,-0.3)(3.95,-0.3)
\rput(3,-0.3){$1$}
\psline[linewidth=0.3pt,linecolor=blue]{<->}(6.3,0)(6.3,0.3)
\rput(6.5,0.15){$\ve$}
\psline[linewidth=0.6pt,linecolor=brown,linestyle=dashed](0,0)(0,2)
}
\end{picture}
\caption{A tube.}
\label{fig tube}
\end{figure}
The lower and upper rows of rectangles are adjacent to the lower and upper sides of the tube, respectively. The distance between neighbor rectangles of each row equals 1. The opening of the tube is the left vertical side of the large rectangle. Denote by $\eta_{\ve,\del,a}$ the measure generated by the tube.

\begin{figure}[h]
\begin{picture}(0,75)
\rput(2.6,0.4){
\scalebox{0.7}{
 \psline[linewidth=1pt,linecolor=brown](0,2.5)(13.51,2.5)
   \psline[linewidth=1pt,linecolor=brown](2.5,2.5)(2.5,2.2)
   \psline[linewidth=1pt,linecolor=brown](5,2.5)(5,2.2)
   \psline[linewidth=1pt,linecolor=brown](7.5,2.5)(7.5,2.2)
   \psline[linewidth=1pt,linecolor=brown](10,2.5)(10,2.2)
   \psline[linewidth=1pt,linecolor=brown](12.5,2.5)(12.5,2.2)
 \psline[linewidth=1pt,linecolor=brown](0,0)(13.5,0)
   \psline[linewidth=1pt,linecolor=brown](2.5,0)(2.5,0.3)
   \psline[linewidth=1pt,linecolor=brown](5,0)(5,0.3)
   \psline[linewidth=1pt,linecolor=brown](7.5,0)(7.5,0.3)
   \psline[linewidth=1pt,linecolor=brown](10,0)(10,0.3)
   \psline[linecolor=brown](12.5,0)(12.5,0.3)
  \psline[linestyle=dotted,linecolor=brown,linewidth=1.1pt](0,0)(0,2.5)
 \psdots[dotsize=1.8pt](13.9,0)(14.3,0)(14.7,0)(15.1,0)(15.5,0)(15.9,0)
 \psdots[dotsize=1.8pt](13.9,2.5)(14.3,2.5)(14.7,2.5)(15.1,2.5)(15.5,2.5)(15.9,2.5)
   \psline[linewidth=1pt,arrows=->,arrowscale=2,linecolor=red](0,0.5)(1.5,1)
   }
}
\end{picture}
\caption{A semi-infinite tube.}
\label{fig semiinfinite tube}
\end{figure}
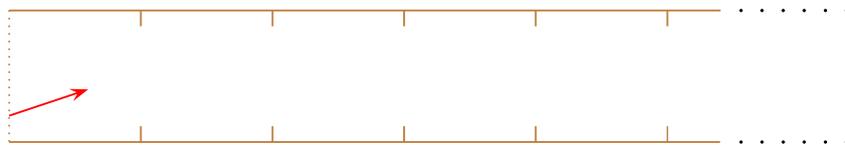

For any particle incident in the tube, with $\vphi$ and $\vphi^+$ being the angles of getting in and getting out, only two cases may happen: $\vphi^+ = \vphi$ or $\vphi^+ = -\vphi$. Letting $a \to \infty$ and $\del \to 0$ (with $\ve$ fixed), we get the semi-infinite tube where small rectangles are substituted with vertical segments of length $\ve$ (see fig. \ref{fig semiinfinite tube}). Studying the dynamics in this tube amounts to the following ergodic problem.

Consider the iterated rotation of the circle by a fixed angle $\al$,\, $\xi_n = \xi + \al n\!\! \mod\! 1$,\, $n = 1,\, 2,\ldots$ and mark the successive moments $n = n_1$,\, $n_1 + n_2$,\, $n_1 + n_2 + n_3,\ldots$, when $\xi_n \in [-\ve,\, \ve]\!\! \mod\! 1$. Denote by $l = l_\ve(\xi,\al)$ the smallest value such that $n_1 - n_2 + \ldots + n_{2l-1} - n_{2l} \le 0$. Let $\mathbb{P}$ be a probability measure on $[0,\, 1] \times [0,\, 1]$ absolutely continuous with respect to Lebesgue measure. Then there exists the limiting distribution $p_k = \lim_{\ve\to 0} \mathbb{P}(\{ (\xi,\al): l_\ve(\xi,\al) = k \})$, with $\sum_{k=1}^\infty p_k = 1$.

In \cite{tube} this statement is proved and is then used to show that that the semi-infinite tube is an asymptotically retroreflecting "hollow"{} (it is not a true hollow, since it is unbounded and its boundary is not piecewise smooth). This means in this case that $\mu(\vphi^+ = \vphi)$ goes to 1
 %$\mu(\{ (\xi,\vphi) \in [0,\, 1] \times \interval : \vphi_\ve^+(\xi,\vphi) = \vphi \}) \to 1$
as $\ve \to 0$.

Let us show that there exists a family of {\it true} tube-shaped hollows which is asymptotically retroreflecting. To this end, define the function $H(\xi,\vphi,\ve,\del,a)$ which is equal to 0, if the billiard particle with the initial data $(\xi,\vphi)$ satisfies the equality $\vphi^+ = \vphi$, and to 1, if $\vphi^+ = -\vphi$ (there are no other possibilities). For the semi-infinite tube this function takes the form $H(\xi,\vphi,\ve,0,+\infty) =: H(\xi,\vphi,\ve)$. The asymptotical retroreflectivity of the semi-infinite tube means that
$$
\lim_{\ve\to 0} \int\!\!\!\int_{[0,\, 1] \times \interval} H(\xi,\vphi,\ve)\, d\xi\, d\vphi = 0.
$$

Note that for fixed $\xi,\ \vphi$ and for $1/a$ and $\del$ small enough the corresponding particle makes the same sequence of reflections (and therefore has the same output velocity) as in the limiting case $\del=0, \ a = +\infty$. This implies that $H(\xi,\vphi,\ve,\del,a)$ pointwise converges (stabilizes) to $H(\xi,\vphi,\ve)$ as $\del \to 0, \ a \to +\infty$, and therefore,
$$
\lim_{\del\to 0,\,a\to+\infty} \int\!\!\!\!\int_{[0,\, 1] \times \interval} \!H(\xi,\vphi,\ve,\del,a)\, d\xi\, d\vphi = \int\!\!\!\!\int_{[0,\, 1] \times \interval} \!H(\xi,\vphi,\ve)\, d\xi\, d\vphi.
$$
Then, using the diagonal method, one selects $\del = \del(\ve)$ and $a = a(\ve)$ such that $\lim_{\ve\to 0} a(\ve) = \infty$,\, $\lim_{\ve\to 0} \del(\ve) = 0$ and
$$
\lim_{\ve\to 0} \int\!\!\!\int_{[0,\, 1] \times \interval} H(\xi,\vphi,\ve,\del(\ve),a(\ve))\, d\xi\, d\vphi = 0.
$$
Thus, the corresponding family of tubes is asymp\-totically retroreflecting.

The obtained result can be formulated as follows.

\begin{theorem}\label{t tube}
$\eta_{\star}$ is a limit point of the set of measures generated by tubes, $\{ \eta_{\ve,\del,a} \}$, equipped with the norm topology.
\end{theorem}

The tube has the following properties.

1. The tube is a {\it convenient} hollow. This property makes it possible to construct an \aprr\ with identical tube-shaped hollows; see fig. \ref{fig Four bodies}(b).

%\begin{figure} \begin{center} \includegraphics[width=0.58\columnwidth]{figTubeBody} \end{center} \caption{A body with tube-shaped hollows.} \label{fig tube body} \end{figure}

2. The measure generated by the tube (with properly chosen $\del = \del(\ve)$ and $a = a(\ve)$) converges {\it in norm} to the \rr\ measure. In other words, the portion of retroreflected particles (that is, particles reflected in the {\it exactly} opposite direction) tends to 1.

3. We believe this construction admits a generalization to higher dimensions, but we could not prove it yet.

4. The average number of reflections in the tube is of the order of $1/\ve$, and therefore, goes to infinity as $\ve \to 0$.

\subsection{\Dopa{}}

Another remarkable hollow called {\it \dopa{}} was discovered and studied by P Gouveia in \cite{Gouv disser} (see also \cite{Gouveia Plakhov Torres AMC}). It is a curvilinear triangle, with the opening being the base of the triangle. Its lateral sides are arcs of parabolas, where the vertex of each parabola coincides with the focus of the other one (and also coincides with a vertex of the triangle at its base). The base is a segment contained in the common axis of the parabolas; see fig. \ref{fig double parabola}.

\begin{figure}[h]
\begin{picture}(0,140)
\rput(2.7,0.1){
\scalebox{1}{
\rput(4.7,0){
\scalebox{6}{
\psecurve[linewidth=0.2pt,linecolor=brown](-0.21875,-0.25)(-0.25,0)(-0.21875,0.25)(-0.125,0.5)(0,0.7071)(0.25,1)
\psecurve[linewidth=0.2pt,linecolor=brown](0.21875,-0.25)(0.25,0)(0.21875,0.25)(0.125,0.5)(0,0.7071)(-0.25,1)
\psline[linewidth=0.2pt,linestyle=dotted,dotsep=0.75pt](-0.25,0)(0.25,0)
}}
   }
}
\end{picture}
\caption{\Dopa{}.}
\label{fig double parabola}
\end{figure}

The \dopa{} is a nearly perfect \rr; the measure $\eta_{\doubpar}$ generated by this hollow satisfies $\FFF(\eta_{\doubpar}) = 0.9977$; this value is only $0.23\%$ smaller than the maximal value of $\FFF$.
A body bounded by \dopa{}s is shown in figure \ref{fig Four bodies}(c).

%\begin{figure} \begin{center} \includegraphics[width=0.44\columnwidth]{rrDoubleParabola} \end{center} \caption{A body with \dopa{}-shaped hollows.} \label{fig double-parabola body} \end{figure}

The \dopa{} has the following properties.

1. It is a {\it convenient} hollow.

\begin{figure}
\begin{center}
\includegraphics[width=0.48\columnwidth]{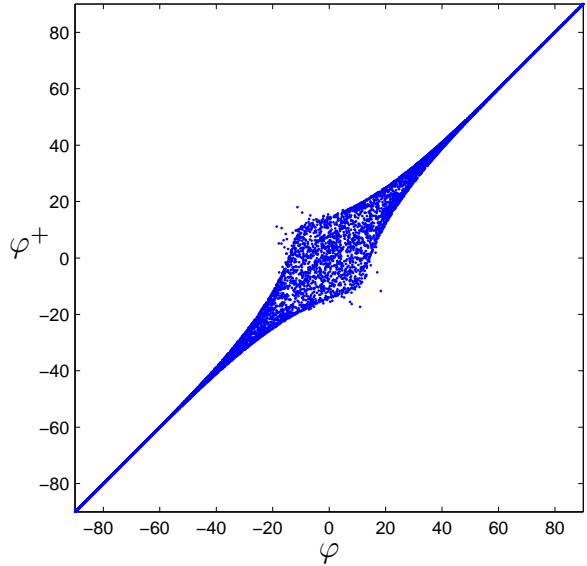}
\end{center}
\caption{The support of the measure generated by \dopa{} is shown. It is obtained numerically by calculating 10\,000 randomly chosen pairs $(\vphi,\vphi^+)$.}
\label{fig distribution phiphi}
 \end{figure}

2. There always exists a small discrepancy between the initial and final directions, which is maximal for perpendicular incidence and vanishes for nearly tangent incidence. See figure \ref{fig distribution phiphi}, where the support of $\eta_{\doubpar}$ is shown. The figure is obtained numerically, by calculating the pairs $(\vphi,\vphi^+)$ for $10\,000$ values of $\vphi$ chosen at random. This means that, when illuminated, the contour of the \rr\ is seen best of all, which is useful for visual reconstruction of its shape.

3. We do not know if there exist multidimensional generalizations of this shape. By now, the greatest value of the parameter $\FFF$ attained by numerical simulation in three dimensions equals $0.9$.

4. For most particles, the number of successive reflections equals 3, although 4, 5, etc. (up to infinity) reflections are also possible. When the number of reflection increases, the number of corresponding particles rapidly decreases.

5. The boundary of \dopa{} is the graph of a function. This means that this shape may be easy for manufacturing.

\subsection{Notched angle}\label{subsec notch ang}

\begin{figure}[h]
\begin{picture}(0,270)
\rput(8,9.4){
\psline[linewidth=0.2pt](-2.4,-9)(0,0)(2.4,-9)
\psline[linewidth=0.2pt](-2,-9)(0,0)(2,-9)
\psline(2,-9)(2,-7.5)(1.66667,-7.5)(1.66667,-6.25)(1.3889,-6.25)(1.3889,-5.20833)(1.15741,-5.20833)
(1.15741,-4.34028)(0.9645,-4.34028)(0.9645,-3.6169)(0.80375,-3.6169)(0.80375,-3.01408)(0.6698,-3.01408)
(0.6698,-2.51174)(0.55816,-2.51174)(0.55816,-2.09311)(0.46514,-2.09311)(0.46514,-1.74426)(0.3876,-1.74426)
(0.3876,-1.45355)(0.323,-1.45355)(0.323,-1.2113)(0.26918,-1.2113)(0.26918,-1.0094)(0.2243,-1.0094)
(0.2243,-0.8412)(0.18693,-0.8412)(0.18693,-0.70098)(0.15577,-0.70098)(0.15577,-0.58415)(0.1298,-0.58415)
(0.1298,-0.4868)(0.10818,-0.4868)
\psline(-2,-9)(-2,-7.5)(-1.6667,-7.5)(-1.6667,-6.25)(-1.3889,-6.25)(-1.3889,-5.2083)(-1.1574,-5.2083)
(-1.15741,-4.34028)(-0.9645,-4.34028)(-0.9645,-3.6169)(-0.80375,-3.6169)(-0.80375,-3.01408)(-0.6698,-3.01408)
(-0.6698,-2.51174)(-0.55816,-2.51174)(-0.55816,-2.09311)(-0.46514,-2.09311)(-0.46514,-1.74426)(-0.3876,-1.74426)
(-0.3876,-1.45355)(-0.323,-1.45355)(-0.323,-1.2113)(-0.26918,-1.2113)(-0.26918,-1.0094)(-0.2243,-1.0094)
(-0.2243,-0.8412)(-0.18693,-0.8412)(-0.18693,-0.70098)(-0.15577,-0.70098)(-0.15577,-0.58415)(-0.1298,-0.58415)
(-0.1298,-0.4868)(-0.10818,-0.4868)
  \psarc(0,0){1.3}{-101.5}{-78.5}
  \rput(0,-1.55){$\al$}
  \psarc(0,0){0.38}{-108.5}{-71.5}
  \psarc(0,0){0.45}{-108.5}{-71.5}
    \rput(0.75,-0.4){$\al+\bt$}
\psline[linewidth=0.4pt,linestyle=dashed](-2,-9)(2,-9)
\rput(-2.05,-9.25){$\AAA$}
\rput(2.05,-9.25){$\BBB$}
\rput(-2.7,-9.25){$\AAA'$}
\rput(2.6,-9.05){$\BBB'$}
\rput(0,0.3){$\OOO$}
  \rput(0,-6){\Huge $\Om$}
\rput(0,-8.7){\large $I$}
\psline[linewidth=0.4pt,linestyle=dashed](-1.6,-2)(1.6,-2)
\rput(-1.7,-2.2){\scalebox{0.8}{$C$}}
\rput(1.7,-2.2){\scalebox{0.8}{$D$}}
}
\end{picture}
\caption{Notched angle.}
\label{fig zub ugol}
\end{figure}
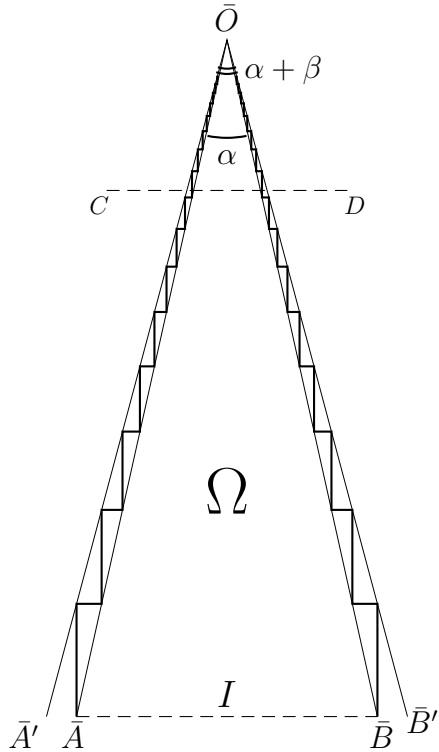

This shape is depicted in figure \ref{fig zub ugol}, and the corresponding body, in figure \ref{fig Four bodies}(d).
Here we point out its properties.

%\vspace{48mm}
%\vspace{12mm}

%\begin{figure} \begin{center} \includegraphics[width=1.10\columnwidth=1.0]{figNotchAngBody} \end{center} \caption{A body with the hollows in the shape of notched angle.} \label{fig body zub ugol} \end{figure}

1. Notched angle is a {\it convenient} hollow.

2. The corresponding measure converges {\it in norm} to the \rr\ measure $\eta_\star$.

3. We are unaware of multidimensional generalizations of this shape.

4. The mean number of reflections in notched angle goes to infinity as $\al$ tends to zero.

5. The boundary of the notched angle is the graph of a function.

The rigorous definition of this shape and the proof of its retroreflectivity are given in the next section \ref{sect notched angle}.

\subsection{Comparison table for \rrs}

Here we put together the billiard \rrs. For convenience, their properties are tabulated below. The limiting values of $r$ are equal to 1 in all shapes, except for the \dopa{}.
\begin{figure}
\begin{center}
\hspace*{-13mm}\includegraphics[width=1.10\columnwidth=1.0]{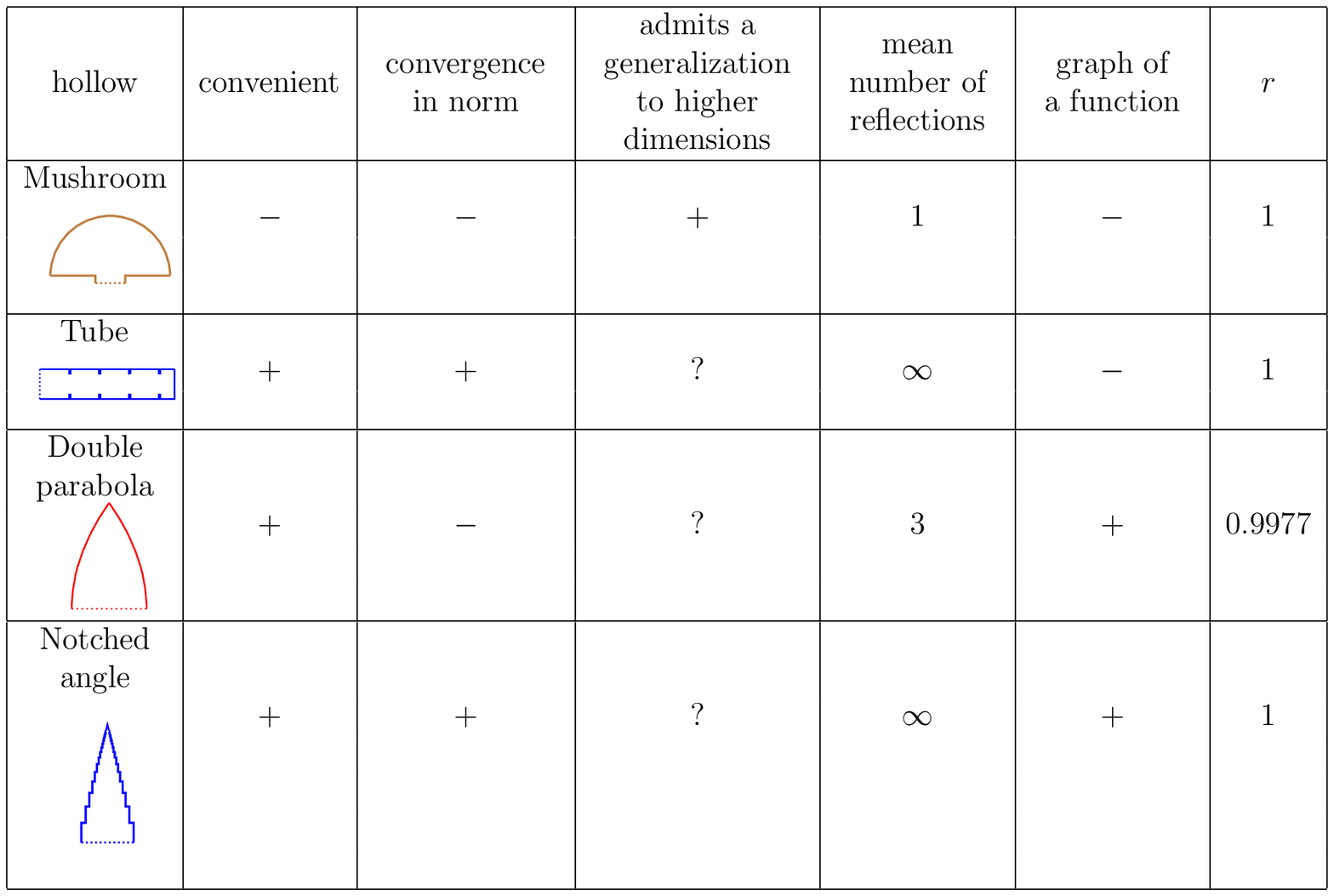}
\end{center}
%\caption{Comparison table for billiard \rrs.}
\label{fig table}
 \end{figure}
In figure \ref{fig Four bodies}, four bodies with boundaries formed by corresponding retroreflecting hollows are shown.

\begin{figure}
\begin{center}
\hspace*{5mm}
\includegraphics[width=1.00\columnwidth]{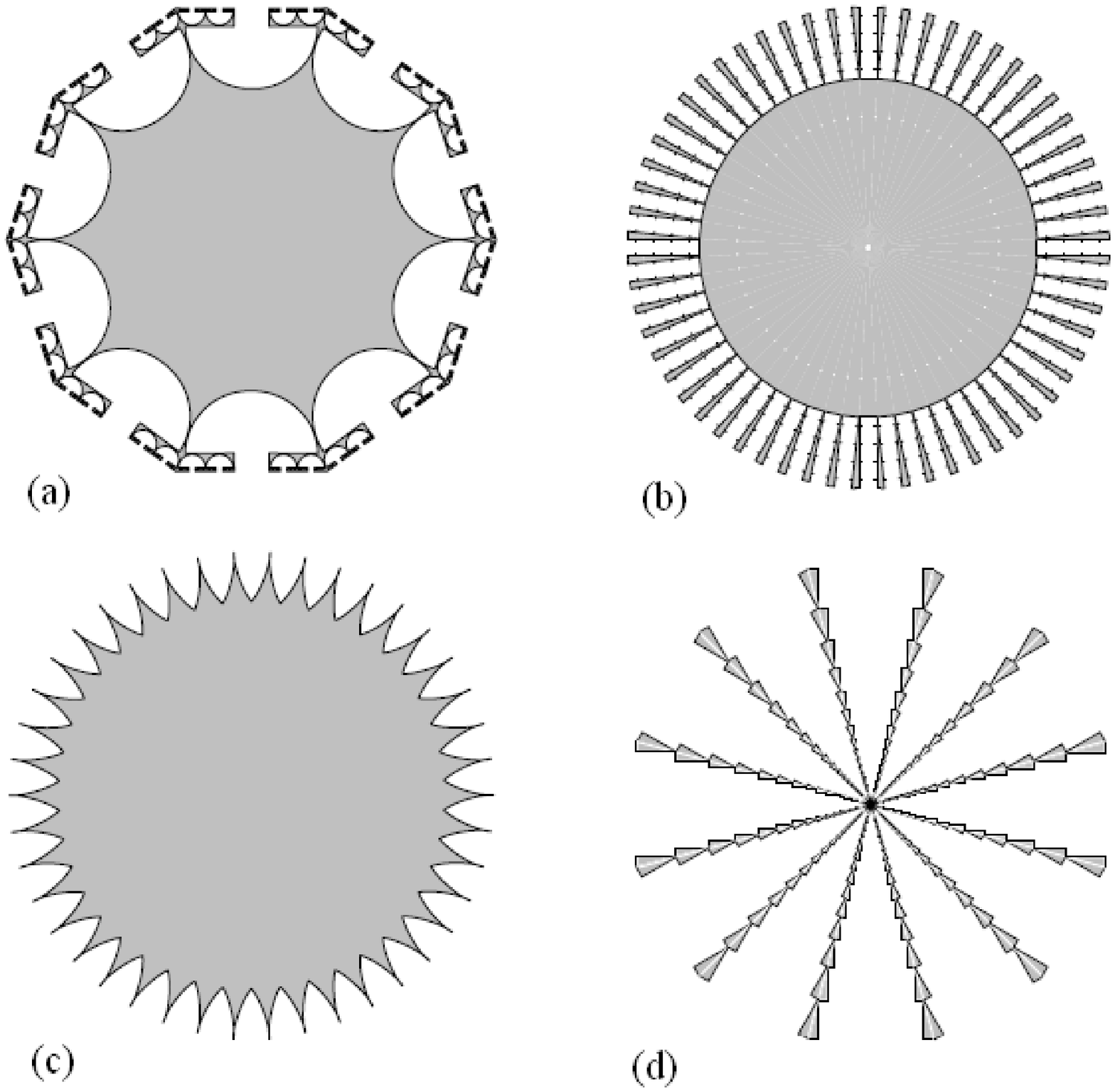}
\end{center}
\caption{Bodies with boundaries formed by retroreflecting hollows: (a) mushroom; (b) tube; (c) \dopa; (d) notched angle.}
\label{fig Four bodies}
\end{figure}

As concerns possible applications of these shapes, each of them seems to have some advantages and disadvantages. Tube and notched angle ensure exact direction reversal, while in mushroom and \dopa\ a small discrepancy between initial and final directions is always present, which can make them inefficient at very large distances. On the other hand, the number of reflections for the most part of particles in mushroom and \dopa\ equals 1 and 3, respectively, while the mean number of reflections goes to infinity for sequences of bodies representing tube and notched angle, which may imply need for high quality of reflecting boundary.

\section{Notched angle}\label{sect notched angle}

Consider two isosceles triangles, $\AAA\OOO\BBB$ and $\AAA'\OOO\BBB'$, with the common vertex $\OOO$ and require that the base of one of them is contained in the base of the other one, $\AAA\BBB \subset \AAA'\BBB'$. The segment $\AAA\BBB$ is horizontal in figure \ref{fig zub ugol}. Denote $\measuredangle \AAA\OOO\BBB = \al$ and $\measuredangle \AAA'\OOO\BBB' = \al + \bt$. %$\al > 0$,\, $\bt > 0$.
Draw two broken lines with horizontal and vertical segments with the origin at $\AAA$ and $\BBB$, respectively, and require that the vertices of the first line belong to the segments $\OOO\AAA$ and $\OOO\AAA'$, and the vertices of the second line, to the segments $\OOO\BBB$ and $\OOO\BBB'$; see figure \ref{fig zub ugol}. The endpoint of both broken lines is $\OOO$; both lines have infinitely many segments and finite length. We will consider the "hollow"{} $(\Om, I)$ with the opening $I = I_\al = \AAA\BBB$ and with the set $\Om = \Om_{\al,\bt}$ bounded by $\AAA\BBB$ and the two broken lines. This "hollow"{} will be called a {\it notched angle} with the size $(\al,\bt)$, or just an $(\al,\bt)$-{\it angle}.
The boundary $\pl\Om$ is not piecewise smooth ($\OOO$ is a limit point for singular points of $\pl\Om$), therefore the word {\it hollow} is put in quotes; however, the measure generated by this "hollow"{} is defined in the standard way. This measure depends only on $\al$ and $\bt$ and is denoted by $\eta_{\al,\bt}$.

\begin{theorem}\label{t zub ugol}
There exists a function $\bt = \bt(\al)$,\, $\lim_{\al\to 0} (\bt/\al) = 0$ such that $\eta_{\al,\bt}$ converges in norm to the \rr\ measure $\eta_{\star}$ as $\al \to 0$.
\end{theorem}

\begin{remark}\label{zam notch true}
{\rm Using this theorem, one easily constructs a family of {\it true} hollows for which convergence in norm to $\eta_{\star}$ takes place. Namely, draw a straight line $CD$ parallel to $\AAA\BBB$ at a small distance $\del$ from $\OOO$; the {\it true hollow} is the part of the original "hollow"{} situated between $\AAA\BBB$ and $CD$, with the same opening (see fig. \ref{fig zub ugol}). The measure generated by this hollow tends to $\eta_{\star}$ as $\al \to 0$, with properly chosen $\bt = \bt(\al)$ and $\del = \del(\al)$ vanishing when $\al \to 0$.}
\end{remark}

\begin{proof}
For any initial data $\xi$,\, $\vphi$ the angle of getting away $\vphi^+ = \vphi^+_{\al,\bt}(\xi,\vphi)$ satisfies either $\vphi^+ = \vphi$, or $\vphi^+ = -\vphi$. To prove the theorem, it suffices to check that the measure $\mu$ of the set of initial data $\xi$,\, $\vphi$ satisfying $\vphi^+_{\al,\bt}(\xi,\vphi) = -\vphi$ and $|\vphi| > \al$ tends to 0 as $\al \to 0$,\, $\bt = \bt(\al)$.

Make a uniform extension along the horizontal axis in such a way that the resulting angle $\AAA\OOO\BBB$ becomes right. Then the angle $\AAA'\OOO\BBB'$ becomes equal to $\pi/2 + \gam$,\, $\gam = \gam(\al,\bt)$ (see fig. \ref{fig red zub ugol}), besides the conditions $\al \to 0$,\, $\bt/\al \to 0$ imply that $\gam \to 0$. This extension takes the $(\al,\bt)$-angle to a $(\pi/2,\del)$-angle, takes each billiard trajectory to another billiard trajectory, and takes the measure $\frac 12 \cos\vphi\,d\vphi\,d\xi$ to a measure absolutely continuous with respect to it.

The vertices of the resulting notched angle will be denoted by $O$,\, $A$,\, $B$,\, $A'$,\, $B'$, without overline, in order to distinguish them from the previous notation.

Without loss of generality we assume that $|OA| = |OB| = 1$. Introduce the uniform parameter $\xi$ on the segment $AB$, where $A$ corresponds to the value $\xi = 0$ and $B$, to the value $\xi = 1$. Extend the trajectory of an incident particle with initial data $\xi$,\, $\vphi < -\pi/4$\footnote{Recall that the angle $\vphi$ is measured counterclockwise from the vertical vector $(0, 1)$ to the velocity of the incident particle, so one has $\vphi < 0$ in figure \ref{fig red zub ugol}.} until the intersection with the extension of  $OA$. Denote by $\tilde x_0$ the distance from $O$ to the point of intersection; see fig. \ref{fig red zub ugol}. (In what follows, a point on the ray $OA$ or $OB$ will be identified with the distance from the vertex $O$ to this point.) In the new representation, the particle starts the motion at a point $\tilde x_0$ and intersects the segment $AB$ at a point $\xi$ and at an angle $\vphi$. Continuing the straight-line motion, it intersects the side $OB$ at a point $x_1$\, $(0 < x_1 < 1)$, then makes one or two reflections from the broken line and intersects $OB$ again at a point $\tilde x_1$. Denote $x_1/\tilde x_0 = \lam$; obviously one has $0 < \lam < 1$. The value $\lam$ is the tangent of the angle of trajectory inclination relative to $OA$; thus, one has $\vphi = -\pi/4 - \arctan\lam$. It is convenient to change the variables in the space of particles getting into the hollow at an angle $\vphi < -\pi/4$. Namely, we pass from the parameters $\xi \in [0,\, 1]$,\, $\vphi \in [-\pi/2,\, -\pi/4]$ to the parameters $\lam \in [0,\, 1]$,\, $\tilde x_0 \in [1,\, 1/\lam]$. This change of variables can be written as $\xi = \frac{\lam}{1 - \lam}\, (\tilde x_0 - 1)$,\, $\vphi = \pi/4 + \arctan\lam$; it transforms the measure $\frac{1}{2}\, \cos\vphi\,d\vphi\,d\xi$ into the measure $\frac{\lam}{2\sqrt{2}(1+\lam^2)^{3/2}}\, d\lam\, d\tilde x_0$.

%\vspace{48mm}
\begin{figure}[h]
\begin{picture}(0,220)
 \scalebox{1.4}{
\rput(5.5,1){
\psline[linewidth=0.14pt](-5,-1)(0,4)(4,0)
\psline[linewidth=0.14pt](-5,0)(0,4)(5,0)
\psline[linewidth=0.6pt,linecolor=blue](4,0)(4,0.8)(3.2,0.8)(3.2,1.44)(2.56,1.44)(2.56,1.952)(2.048,1.952)(2.048,2.3616)(1.6384,2.3616)
(1.6384,2.68928)(1.31072,2.68928)(1.31072,2.951424)(1.048576,2.951424)(1.048576,3.1611392)(0.8388608,3.1611392)
(0.8388608,3.32891136)(0.67108864,3.32891136)(0.67108864,3.463129088)(0.536870912,3.463129088)(0.536870912,3.57050327)
(0.429497,3.570503)(0.429497,3.656403)(0.343597,3.656403)
(0.343597,3.725122)(0.274878,3.725122)(0.274878,3.780098)
(0.219902,3.780098)(0.219902,3.824078)
(0.175922,3.824078)(0.175922,3.859262)
\psline[linewidth=0.6pt,linecolor=blue](-4,0)(-4,0.8)(-3.2,0.8)(-3.2,1.44)(-2.56,1.44)(-2.56,1.952)(-2.048,1.952)(-2.048,2.3616)(-1.6384,2.3616)
(-1.6384,2.68928)(-1.31072,2.68928)(-1.31072,2.951424)(-1.048576,2.951424)(-1.048576,3.1611392)(-0.8388608,3.1611392)
(-0.8388608,3.32891136)(-0.67108864,3.32891136)(-0.67108864,3.463129088)(-0.536870912,3.463129088)(-0.536870912,3.57050327)
(-0.429497,3.570503)(-0.429497,3.656403)(-0.343597,3.656403)
(-0.343597,3.725122)(-0.274878,3.725122)(-0.274878,3.780098)
(-0.219902,3.780098)(-0.219902,3.824078)
(-0.175922,3.824078)(-0.175922,3.859262)
%\psline[linewidth=0.4pt](-1.7,3.57050327)(1.7,3.57050327)
\psline[linewidth=0.3pt,linestyle=dashed,linecolor=blue](-4,0)(4,0)
\rput(-4.2,0.07){\scalebox{0.7}{$A$}}
\rput(4.12,-0.17){\scalebox{0.7}{$B$}}
\rput(-5.2,-0.15){\scalebox{0.7}{$A'$}}
\rput(5.15,-0.12){\scalebox{0.7}{$B'$}}
\rput(0,4.2){\scalebox{0.7}{$O$}}
 \rput(-4.55,-0.8){\scalebox{0.7}{$\tilde x_0$}}
 \rput(2.87,0.86){\scalebox{0.7}{$x_1$}}
  \rput(-1.7,2.03){\scalebox{0.7}{$x_2$}}
    \rput(2.65,1.13){\scalebox{0.6}{$\tilde x_1$}}
        \rput(-1.4,2.4){\scalebox{0.6}{$\tilde x_2$}}
 \psdots[dotsize=1.4pt](-4.58,-0.58)(2.93,1.07)(2.77,1.23)(-1.77,2.23)(-1.66,2.34)(-1.9,0)
 \rput(-1.9,-0.2){\scalebox{0.7}{$\xi$}}
 \psline[linewidth=0.4pt,linecolor=brown](-1.9,0)(-1.9,0.4)
 \psarc[linewidth=0.4pt,linecolor=brown](-1.9,0){0.2}{12.5}{90}
 \rput(-1.7,0.3){\scalebox{0.7}{$\vphi$}}
\psline[linewidth=0.3pt,linecolor=red,arrows=->,arrowscale=1.5](-4.6,-0.6)(-0.7,0.265)
\psline[linewidth=0.3pt,linecolor=red,arrows=->,arrowscale=1.5](1.4507,1.518)(-2.048,2.294)(-1.75,2.3616)(0.2,1.9291)
\psline[linewidth=0.3pt,linecolor=red,arrows=->,arrowscale=1.5](-0.7,0.265)(3.2,1.13)(1.4507,1.518)
      %\rput(2.4,1.3){$\tilde x_1$}
%\rput(-2,3.57){$C$}
%\rput(2,3.57){$D$}
  \psarc[linewidth=0.9pt](0,4){0.99}{218}{226}
  \psarc[linewidth=0.9pt](0,4){0.92}{218}{226}
  \rput(-0.9,3.56){\scalebox{0.65}{$\gam/2$}}
}
 }
\end{picture}
\caption{The reduced notched angle.}
\label{fig red zub ugol}
\end{figure}
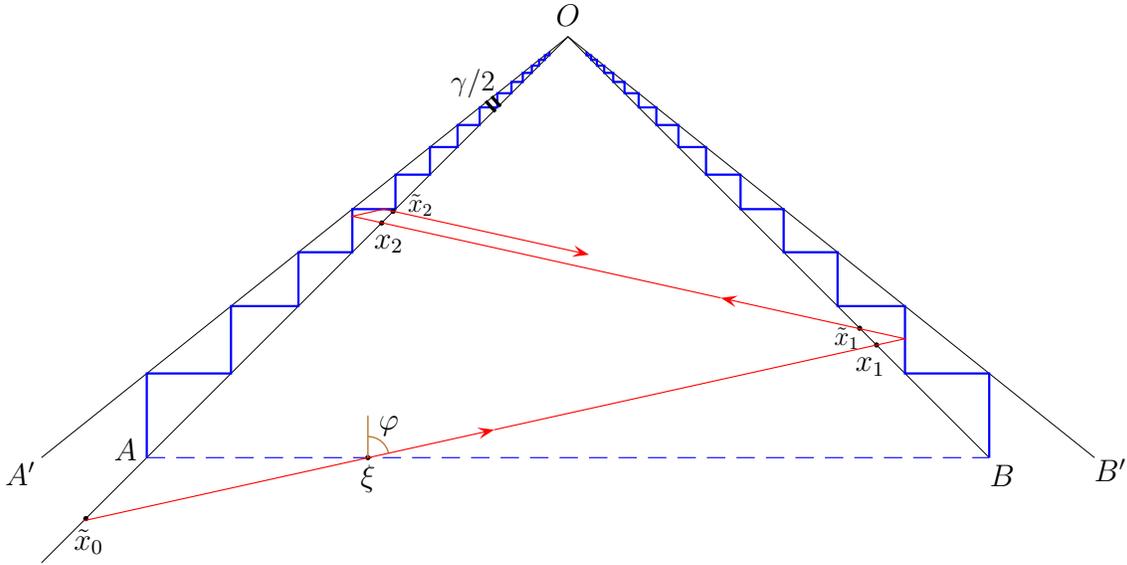
%\vspace{12mm}

By considering successive alternating reflections of the particle from the broken lines resting on the sides $OB$ and $OA$, we define the sequence of values $x_1$,\, $\tilde x_1$, \ldots, $x_{m-1}$,\, $\tilde x_{m-1}$. Obviously, all these values are smaller than 1. Then the particle gets out of the hollow and intersects the extension of the side $OA$ or $OB$ at a point $x_m > 1$. If $m$ is even, then the intersection with $OA$ takes place, and $\vphi^+ = \vphi$. If $m$ is odd, then intersection with $OB$ takes place, with  $\vphi^+ = -\vphi$. Clearly, $m$ depends on the initial data $\tilde x_0$,\, $\lam$ and on the parameter $\gam$,\, $m = m_\gam(\tilde x_0, \lam)$.

\begin{utv}\label{utv 1}
For any $\lam$, the measure of the set of values $\tilde x_0$ such that $m_\gam(\tilde x_0, \lam)$ is odd, goes to $0$ as $\gam \to 0$.
\end{utv}

Let us derive the theorem from this statement. Indeed, let $f_\gam(\lam)$ be the measure of the set indicated in the statement, $f_\gam(\lam) = |\{ \tilde x_0 : m_\gam(\tilde x_0, \lam) \text{ is odd}\, \}|$. Introduce the measure $\eta$ on the segment $[0,\, 1]$ according to $d\eta(\lam) = \frac{\lam\, d\lam}{2\sqrt{2}(1 + \lam^2)^{3/2}}$; then $\int_0^1 f_\gam(\lam)\, d\eta(\lam)$ is the measure of the set of initial values $(\lam,\tilde x_0)$ such that $m_\gam(\tilde x_0, \lam)$ is odd. The value $f_\gam(\lam)$ does not exceed the full Lebesgue measure of the segment $[1,\, \lam^{-1}]$,
 \beq\label{ineq 1}
f_\gam(\lam) \le \lam^{-1} - 1,
 \eeq
and the function $\lam^{-1} - 1$ is integrable relative to $\eta$, $\int_0^1 (\lam^{-1} - 1)\, d\eta(\lam) = \frac{\sqrt{2} - 1}{2\sqrt{2}}$. According to statement \ref{utv 1}, for any $\lam$ holds
 \beq\label{limit 2}
\lim_{\gam\to 0} f_\gam(\lam) = 0.
 \eeq
Taking into account (\ref{ineq 1}) and (\ref{limit 2}) and applying Lebesgue's dominated convergence theorem, one gets
 \begin{equation*}\label{limit 3}
\lim_{\gam\to 0} \int_0^1 f_\gam(\lam)\, d\eta(\lam) = 0.
 \end{equation*}
This means that the measure of the set of values $(\xi, \vphi)$,\, $\vphi \le -\pi/4$ for which the equality $\vphi^+_{\pi/2,\gam}(\xi,\vphi) = -\vphi$ is valid, tends to 0 as $\gam \to 0$. The same statement, due to the axial symmetry of the billiard, is also valid for $\vphi \ge \pi/4$.

Now make a uniform contraction along the abscissa axis transforming the $(\pi/2,\gam)$-angle into an $(\al,\bt)$-angle (where $\bt$ depends on $\gam$ and $\al$). Taking into account that the measures generated by these angles are mutually absolutely continuous, we get that the measure $\mu(\{ (\xi, \vphi) : |\vphi| \ge \al \text{ and } \vphi^+_{\al,\bt} (\xi,\vphi) = -\vphi \})$ goes to 0 at fixed $\al$ and $\bt \to 0$.

Finally, choose a diagonal family of parameters $\al$,\, $\bt(\al)$,\, $\lim_{\al\to 0}(\bt(\al)/\al) = 0$ such that the measure
$$
\mu(\{ (\xi, \vphi) : |\vphi| \ge \al \, \text{ and } \, \vphi^+_{\al,\bt(\al)}(\xi,\vphi) = -\vphi \}) \to 0
\quad \text{as} \ \ \al \to 0.
$$
It remains to notice that $\mu(\vphi^+_{\al,\bt(\al)} = -\vphi) \le \mu(|\vphi| \ge \al \text{ and } \vphi^+_{\al,\bt(\al)} = -\vphi) + \mu(|\vphi| < \al)$ and $\mu(|\vphi| < \al) \to 0$ as $\al \to 0$. This finishes the proof of theorem \ref{t zub ugol}.
\end{proof}

{\it Proof of statement \ref{utv 1}}.
Note that the broken lines intersect with the sides $OA$ and $OB$ at the points $x = e^{-n\del}$,\, $n = 0,\, 1,\, 2, \ldots$, where $\del$ is defined by the relation $\tanh\del = \sin\gam$. Consider an arbitrary pair of values $x_k$,\, $\tilde x_k$; they belong to a segment bounded by a pair of points $x = e^{-n\del}$ and $e^{-(n+1)\del}$. Consider also the right triangle, with the hypotenuse being this segment and with the legs being segments of the broken line.

\begin{figure}[h]
\begin{picture}(0,120)
\rput(7.8,2.2){
\pspolygon[linewidth=0.2pt](-2,0)(0,2)(2,0)(0,-2)
\psline[linewidth=0.8pt](2,0)(0,0)(0,-2)
\psline[linewidth=0.2pt](-2,0)(0,0)(0,2)
\psline[linewidth=0.4pt](2,0)(0,-2)
\rput(0.55,-2.2){$\zeta = 0$}
\rput(2.58,0){$\zeta = 1$}
\psline[linewidth=0.6pt,linestyle=dashed](-2,0)(1.25,-0.75)
\psdots[dotsize=2.5pt](2,0)(0,-2)(1.25,-0.75)(-1.8,0.2)(-1.52,-0.48)
\rput(1.85,-0.85){$\zeta = \lam$}
\psline[linewidth=0.4pt](-1.8,0.2)(-0.76,-0.04)
\psline[linewidth=0.4pt,arrows=<-,arrowscale=1.5](-0.76,-0.04)(1.45,-0.55)
\psline[linewidth=0.4pt](-1.52,-0.48)(-1,-0.6)
\psline[linewidth=0.4pt,arrows=<-,arrowscale=1.5](-1,-0.6)(0.95,-1.05)
\rput(-2.7,0.45){$1 + \lam - \zeta$}
\rput(-2.1,-0.7){$\lam^{-1} \zeta$}
}
\end{picture}
\caption{Dynamics in a small right triangle.}
\label{fig auxiliary}
\end{figure}
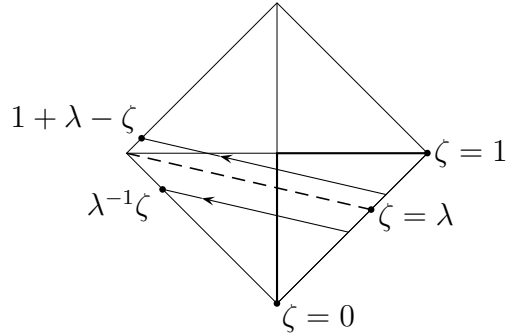

Two cases may happen: either (I) $x_k/\tilde x_{k-1} = \lam$ or (II) $x_k/\tilde x_{k-1} = \lam^{-1}$, the first case corresponding to the "forward"{} motion in the direction of the point $O$, and the second, to the "backward"{} motion. Introduce the local variable $\zeta$ on the hypotenuse according to $x = e^{-n\del} [1 + \zeta(e^{-\del} - 1)]$ (see fig. \ref{fig auxiliary}). Thus, the value $\zeta = 0$ corresponds to the point $x = e^{-n\del}$, and $\zeta = 1$, to the point $x = e^{-(n+1)\del}$. The sequences $x_k$,\, $\tilde x_k$ generate two sequences $\zeta_k$,\, $\tilde \zeta_k \in (0,\, 1)$ and an integer-valued sequence $n_k$. Consider the two cases separately.
 \vspace{1mm}
 %\begin{quotation}

%\hspace*{-6.5mm}{\bf 1.}
(I) $x_k/\tilde x_{k-1} = \lam$.

\hspace*{5mm}(a) If $0 < \zeta_k < \lam$, then $\tilde\zeta_k = \lam^{-1} \zeta_k$ and the particle, after leaving the triangle, continues the forward motion, that is, $x_{k+1}/\tilde x_{k} = \lam$.

\hspace*{5mm}(b) If $\lam < \zeta_k < 1$, then $\tilde\zeta_k = 1 + \lam - \zeta_k$ and the particle, after leaving the triangle, proceeds to the backward motion, $x_{k+1}/\tilde x_{k} = \lam^{-1}$.
 \vspace{1mm}

%\hspace*{-6.5mm}{\bf 2.}
(II) $x_k/\tilde x_{k-1} = \lam^{-1}$. In this case one has $\tilde\zeta_k = \lam \zeta_k$ and the backward motion continues, $x_{k+1}/\tilde x_{k} = \lam^{-1}$.
  \vspace{1mm}
     %\end{quotation}

Introduce the logarithmic scale $z = -\frac{1}{\del}\, \ln x$; then one gets a sequence of values $\tilde z_0$,\, $z_1$,\, $\tilde z_1, \ldots, z_{m-1}$,\, $\tilde z_{m-1}$,\, $z_m$. The first and the last term in this sequence are negative, and the rest of the terms are positive. One has $-\frac{1}{\del} \ln \frac{1}{\lam} < \tilde z_0 < 0$. The following equations establish the connection between $z_k, \ \tilde z_k$ and $\zeta_k, \ \tilde \zeta_k$.
 \beq\label{dynamics 1}
z_k = n_k - \frac{1}{\del} \ln[1 + \zeta_k(e^{-\del} - 1)],
 \eeq
 \beq\label{dynamics 2}
\tilde z_k = n_k - \frac{1}{\del} \ln[1 + \tilde \zeta_k(e^{-\del} - 1)].
 \eeq
As $\del \to 0$, one gets $z_k = n_k + \zeta_k + O(\del)$,\, $\tilde z_k = n_k + \tilde \zeta_k + O(\del)$, where the estimates $O(\del)$ are uniform over all $k$ and all initial data; thus, $\zeta_k$ and $\tilde \zeta_k$ are approximately equal to the fractional parts of $z_k$ and $\tilde z_k$, respectively.

For several initial values $k = 1,\, 2, \ldots, k_\del - 1$ corresponding to the forward motion of the particle, according to (Ia) one has
 \beq\label{dynamics 3}
z_k = \tilde z_{k-1} + \frac{1}{\del} \ln \frac{1}{\lam}; \quad \quad \quad 0 < \zeta_k < \lam, \quad \tilde\zeta_k = \lam^{-1} \zeta_k; \quad \quad \quad z_{k+1} = \tilde z_{k} + \frac{1}{\del} \ln \frac{1}{\lam}.
 \eeq
Here and in the following formulas (\ref{dynamics 4}),(\ref{dynamics 5}), $\,\zeta_k$ is determined by $z_k$ and $\tilde z_k$ is determined by $\tilde \zeta_k$, according to (\ref{dynamics 1}) and (\ref{dynamics 2}).
For the value $k = k_\del$ corresponding to the transition from the forward motion to the backward one, according to (Ib) one has
 \beq\label{dynamics 4}
z_{k_\del} = \tilde z_{k_\del-1} + \frac{1}{\del} \ln \frac{1}{\lam}; \ \quad \lam < \zeta_{k_\del} < 1, \quad \tilde\zeta_{k_\del} = 1 + \lam - \zeta_{k_\del}; \quad \ z_{k_\del + 1} = \tilde z_{k_\del} - \frac{1}{\del} \ln \frac{1}{\lam}.
 \eeq
Finally, for the values $k = k_\del + 1, \ldots, m - 1$ corresponding to the backward motion, according to (II) one has
 \beq\label{dynamics 5}
z_{k} = \tilde z_{k-1} - \frac{1}{\del} \ln \frac{1}{\lam}; \quad \quad \quad \tilde\zeta_k = \lam \zeta_k; \quad \quad \quad z_{k+1} = \tilde z_{k} - \frac{1}{\del} \ln \frac{1}{\lam}.
 \eeq
Notice that in figure \ref{fig red zub ugol} one has $k_\del = 2$.

The formulas (\ref{dynamics 1})--(\ref{dynamics 5}) define iterations of the pairs of mappings
\beq\label{pair of mappings}
\tilde z_{k-1} \mapsto z_k \mapsto \tilde z_k
\eeq
with positive integer time $k$. These mappings commute with the shift $z \mapsto z+1$. The initial value $\tilde z_0$ satisfies $\tilde z_0 \in (-\frac{1}{\del} \ln \frac{1}{\lam},\, 0)$, and the relation $z_m \in (-\frac{1}{\del} \ln \frac{1}{\lam},\, 0)$ defines the time $m$ when the corresponding value leaves the positive semi-axis $z \ge 0$ and the process stops.\footnote{Notice that $m$ depends on $\del$ and $\tilde z_0$; thus, strictly speaking, one should write $m = m_\del(\tilde z_0)$. Then the equality holds $m_\del(\tilde z_0) = m_\gam(\tilde x_0, \lam)$, where $\sin\gam = \tanh\del$ and $\tilde x_0 = e^{-\del\tilde z_0}$; recall that the parameter $\lam$ is fixed.}

During the forward motion, the first mapping in (\ref{pair of mappings}) increases the value of $z$ by $\frac{1}{\del} \ln \frac{1}{\lam}$, and the second one changes it by a value smaller than 1. During the backward motion, the first mapping decreases $z$ by $\frac{1}{\del} \ln \frac{1}{\lam}$, and the second mapping changes it again by a value smaller than 1. Therefore, if the initial value satisfies $\tilde z_0 \in (-\frac{1}{\del} \ln \frac{1}{\lam} + 2k,\, -2k)$ with $k > k_\del$, then $z_{2k_\del} \in (-\frac{1}{\del} \ln \frac{1}{\lam},\, 0)$, and so, $m = 2k_\del$. This means that $m$ is always even, except for a small portion $4k/(\frac{1}{\del}\ln\frac{1}{\lam})$ of the initial values. Thus, to complete the proof of statement \ref{utv 1}, we only need a result stating that the transition time $k_\del$ remains bounded when $\del \to 0$.

Due to invariance with respect to integer shifts, the formulas (\ref{dynamics 1})--(\ref{dynamics 5}) determine iterated maps on the unit circumference with the coordinate $z \!\!\mod \!1$. The value $k_\del = k_\del(\tilde z_0 \!\!\mod \!1)$ is a Borel measurable function; it can be interpreted as a random variable, where the random event is represented by the variable $\tilde z_0 \!\!\mod \!1$ on the circumference with Lebesgue measure.

\begin{utv}\label{utv 2}
The limiting distribution of $k_\del$ as $\del \to 0$ equals $P_\lam(k) = \lam^{k-1} (1 - \lam)$,\, $k = 1,\, 2, \ldots$.
\end{utv}

Let us derive statement \ref{utv 1} using statement \ref{utv 2}. Indeed, one has $1 - P_\lam(1) - \ldots - P_\lam(k) = \lam^k$. Take an arbitrary $\ve > 0$ and choose $k$ such that $\lam^k < \ve$. Then, using statement \ref{utv 2}, choose $\del_0 > 0$ such that $\mathbb{P}(k_\del > k) < \ve$ for any $\del < \del_0$. This implies that the inequality $|\tilde z_0 - z_{2k_\del}| < 2k$ holds with the probability at least $1 - \ve$. Therefore, if $\del$ satisfies $\del < \del_0$ and $4k/(\frac{1}{\del}\ln\frac{1}{\lam}) < \ve$, the relative Lebesgue measure of the set of points $\tilde z_0 \in (-\frac{1}{\del}\ln \frac{1}{\lam},\, 0)$ producing the value $m = 2k_\del$ is greater than $1 - 2\ve$. Passing from the variable $\tilde z_0$ to the variable $\tilde x_0 = e^{-\del\tilde z_0}$, one concludes that Lebesgue measure of the set of values of $\tilde x_0$ corresponding to odd $m$ tends to 0 as $\del \to 0$. This completes the proof of statement \ref{utv 1}.~~$\Box$
 \vspace{1mm}

{\it Proof of statement \ref{utv 2}}.
For convenience write down iterations of the pair of mappings until the transition time $k_\del$ in the form
 \beq\label{iterations}
z_k = \tilde z_{k-1} + \frac{1}{\del}\ln\frac{1}{\lam}\!\! \mod\! 1, \quad \tilde z_{k} =f_\del^{-1}(z_k) \quad (1 \le k < k_\del),
 \eeq
where the function $f_\del$ is given by relations (\ref{dynamics 1}), (\ref{dynamics 2}) and (\ref{dynamics 3}); one easily derives that $f_\del(\tilde z) = \zeta^{-1}(\lam\,\zeta(\tilde z))$, with $\zeta(z) = (1 - e^{-\del z})/(1 - e^{-\del})$. The function $f_\del$ is monotone and injectively maps the circumference $\RRR/\ZZZ$ with the coordinate $z\!\! \mod\! 1$ into itself, and is discontinuous at $0\!\! \mod\! 1$. In the limit $\del \to 0$,\, $f_\del(\tilde z)$ uniformly converges to $\lam\tilde z$ and the derivative $f'_\del$ uniformly converges to $\lam$; the last means that
\beq\label{uniform}
\lim_{\del\to 0} \inf f'_\del = \lim_{\del\to 0} \sup f'_\del = \lam.
\eeq
The iterations (\ref{iterations}) are defined while $z_k \in \text{Range}(f_\del)$; the first moment when $z_k \not\in \text{Range}(f_\del)$ is $k = k_\del$.

Denote by $\mathcal{A}_\del(k) = \{ \tilde z_0\!\! \mod\! 1 :\, k_\del(\tilde z_0\!\! \mod\! 1) > k \}$ the set of initial values $\tilde z_0\!\! \mod\! 1 \in \RRR/\ZZZ$ for which the inequality $k_\del > k$ holds true. Then one has $\mathbb{P}(k_\del > k) = |\mathcal{A}_\del(k)|$, where $|\cdot|$ means Lebesgue measure on $\RRR/\ZZZ$. The following inductive formulas are valid: $\mathcal{A}_\del(0) = \RRR/\ZZZ$ and $\mathcal{A}_\del(k+1) = f_\del(\mathcal{A}_\del(k)) - \frac{1}{\del}\ln\frac{1}{\lam}\!\! \mod\! 1$. They imply that $|\mathcal{A}_\del(0)| = 1$ and
 \beq\label{compress}
\inf_z f'_\del(z) \le \frac{|\mathcal{A}_\del(k+1)|}{|\mathcal{A}_\del(k)|} \le \sup_z f'_\del(z).
 \eeq
Formulas (\ref{uniform}) and (\ref{compress}) imply that $\lim_{\del\to 0} |\mathcal{A}_\del(k)| = \lam^k$;
therefore $\lim_{\del\to 0} \mathbb{P}(k_\del = k) = \lim_{\del\to 0} (\mathcal{A}_\del(k-1)) - \mathcal{A}_\del(k))) = \lam^{k-1} (1 - \lam)$. Statement \ref{utv 2} is proved.~~~~~ $\Box$

\section{Appendix}

\subsection{Convergence of measures generated by rectangular hollows}\label{app 1}

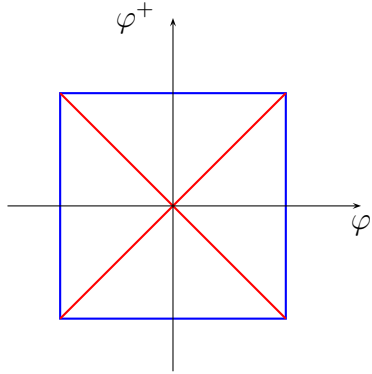
\begin{figure}[h]
\begin{picture}(0,130)
\rput(8,2){
\pspolygon[linecolor=blue](-1.5,-1.5)(1.5,-1.5)(1.5,1.5)(-1.5,1.5)
\psline[linecolor=red](1.5,-1.5)(-1.5,1.5)
\psline[linecolor=red](1.5,1.5)(-1.5,-1.5)
\psline[linewidth=0.4pt]{->}(-2.2,0)(2.5,0)
\psline[linewidth=0.4pt]{->}(0,-2.2)(0,2.5)
\rput(2.5,-0.25){$\vphi$}
\rput(-0.5,2.5){$\vphi^+$}
}
\end{picture}
\caption{The support of the semi-retroreflecting measure.}
\label{fig cross-shaped support}
\end{figure}

Both the measures $\eta^\ve_\sqcup$ and the limiting measure $\frac 12 (\eta_0 + \eta_\star)$  have a cross-shaped support, as shown in figure \ref{fig cross-shaped support}.
Therefore, the density of $\eta^\ve_\sqcup$ can be written down as
$$
\rho_\ve(\vphi)\, \del(\vphi-\vphi^+) + \Big( \frac 12 \cos\vphi - \rho_\ve(\vphi) \Big)\, \del(\vphi+\vphi^+),
$$
and the density of $\frac 12 (\eta_0 + \eta_\star)$ equals
$$
\frac 14 \cos\vphi\, (\del(\vphi-\vphi^+) + \del(\vphi+\vphi^+)).
$$

\begin{figure}[h]
\begin{picture}(0,120)
\rput(5,0.5){
\pspolygon[linecolor=white,fillstyle=solid,fillcolor=yellow](-1,0)(3,0)(3,1)(-1,1)
\psline[linecolor=brown](-1,0)(7,0)
\psline[linecolor=brown](-1,1)(7,1)
\psline[linecolor=brown](-1,2)(7,2)
\psline[linecolor=brown](-1,3)(7,3)
\psline[linecolor=brown](3,0)(3,3)
\psline[linecolor=brown,linestyle=dashed,linewidth=0.4pt](-1,0)(-1,3)
\psline[linecolor=brown,linestyle=dashed,linewidth=0.4pt](7,0)(7,3)
\psline[linecolor=red,arrows=->,arrowscale=2](1.5,0.95)(7,2.6)
\psline[linecolor=red,arrows=->,arrowscale=2](-1,0.2)(1.5,0.95)
\psline[linecolor=red,linewidth=0.4pt,linestyle=dashed,arrows=->,arrowscale=2]
(1.667,1)(3,0.6)(1,0)
\psline[linecolor=red,linewidth=0.4pt,linestyle=solid,arrows=->,arrowscale=2]
(1,0)(-1,0.6)
%\rput(2.5,-0.25){$\vphi$}
%\rput(-0.5,2.5){$\vphi^+$}
}
\end{picture}
\caption{The unfolded billiard trajectory in the tube.}
\label{fig unfold in rectangle}
\end{figure}
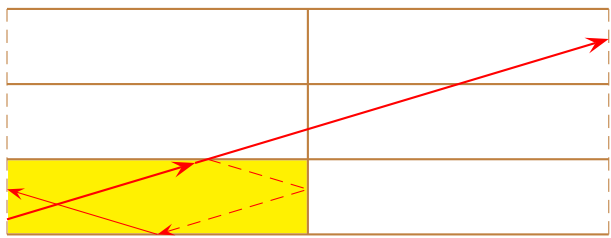

Define the function $f_\ve(\xi,\vphi) = \left\{ \begin{array}{rl} 1, & \text{if } \ \vphi^+(\xi,\vphi) = \vphi\\
-1, & \text{if } \ \vphi^+(\xi,\vphi) = -\vphi \end{array} \right.$;
then one has
$$
\rho_\ve(\vphi) - \Big( \frac 12 \cos\vphi - \rho_\ve(\vphi) \Big) = \cos\vphi \cdot \int_0^1 f_\ve(\xi,\vphi)\, d\xi.
$$
The value of $f_\ve$ is determined from the parity of the number of reflections in the tube and can be easily found by unfolding of the billiard trajectory (see fig. \ref{fig unfold in rectangle}).
One easily sees that $f(\xi,\vphi) = 1$, if $\lfloor \xi + \frac 2\ve \tan\vphi \rfloor$ is odd and $f(\xi,\vphi) = -1$, if $\lfloor \xi + \frac 2\ve \tan\vphi \rfloor$ is even, where $\lfloor \ldots \rfloor$ means the integer part of a real number.

To prove the weak convergence, it suffices to check that for any $-\pi/2 < \Phi_1 < \Phi_2 < \pi/2$,
\beq\label{appendix convergence}
\lim_{\ve\to 0} \int_0^1\!\!\!\int_{\Phi_1}^{\Phi_2} f_\ve(\xi,\vphi)\, \cos\vphi\, d\vphi\, d\xi = 0.
\eeq

Fix $\xi$ and denote $\vphi_m = \arctan(\frac\ve 2(m-\xi))$. One has $f_\ve(\xi,\vphi) = 1$, if $\vphi_{2n-1} < \vphi < \vphi_{2n}$ and $f_\ve(\xi,\vphi) = -1$, if $\vphi_{2n} < \vphi < \vphi_{2n+1}$. One easily deduces from this that the integral $\int_{\Phi_1}^{\Phi_2} f_\ve(\xi,\vphi)\, \cos\vphi\, d\vphi$ converges to zero as $\ve \to 0$ (and is obviously bounded, $|\int_{\Phi_1}^{\Phi_2} f_\ve(\xi,\vphi)\, \cos\vphi\, d\vphi| < 2$), and therefore, the convergence in (\ref{appendix convergence}) takes place.

\subsection{Convergence of measures generated by triangular hollows}\label{app 2}

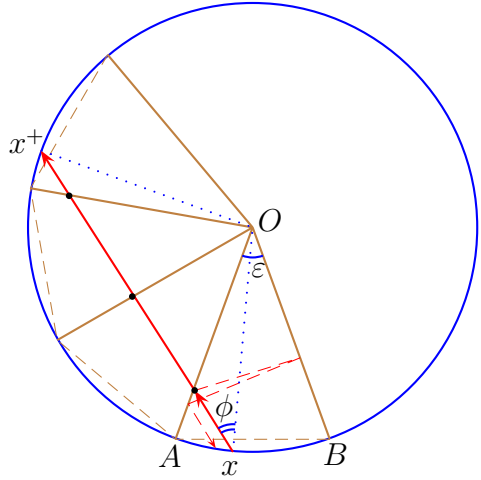
\begin{figure}[h]
\begin{picture}(0,190)
\rput(8.1,3.5){
\pscircle[linecolor=blue](0,0){3}
\psline[linecolor=brown](0,0)(1.026,-2.819)
\psline[linecolor=brown](0,0)(-1.026,-2.819)
\psline[linecolor=brown](0,0)(-2.598,-1.5)
\psline[linecolor=brown](0,0)(-2.954,0.521)
\psline[linecolor=brown](0,0)(-1.928,2.298)
\psline[linecolor=brown,linestyle=dashed,linewidth=0.4pt]
(1.026,-2.819)(-1.026,-2.819)(-2.598,-1.5)(-2.954,0.521)(-1.928,2.298)
\psline[linecolor=red,arrows=->,arrowscale=1.5](-0.77,-2.17)(-2.819,1.026)
\psline[linecolor=red,arrows=->,arrowscale=1.5](-0.261,-2.988)(-0.77,-2.17)
\psline[linecolor=red,linewidth=0.2pt,linestyle=dashed,arrows=->,arrowscale=1.5]
(-0.77,-2.17)(0.631,-1.735)(-0.855,-2.349)(-0.49,-2.94)
\psline[linecolor=blue,linestyle=dotted,linewidth=0.8pt](-2.819,1.026)(0,0)(-0.261,-2.988)
\rput(-0.3,-3.2){$x$}
\rput(-3,1.15){$x^+$}
\psarc[linecolor=blue](-0.261,-2.988){0.3}{85}{122}
\psarc[linecolor=blue](-0.261,-2.988){0.37}{85}{122}
\rput(-0.38,-2.4){$\phi$}
\rput(0.23,0.1){$O$}
\rput(-1.1,-3.05){$A$}
\rput(1.1,-3.02){$B$}
\psdots[dotsize=2.5pt](-0.77,-2.17)(-1.6,-0.92)(-2.44,0.42)
\psarc[linecolor=blue](0,0){0.4}{250}{290}
\rput(0.08,-0.57){$\ve$}
}
\end{picture}
\caption{The unfolded billiard trajectory in the triangle.}
\label{fig unfold triangle}
\end{figure}

The images of the triangular hollow $AOB$ obtained by the unfolding procedure form a polygon inscribed in a circle (see figure \ref{fig unfold triangle}). Introduce the angular coordinate $x \!\! \mod \! 2\pi$ (measured clockwise from the point $B$) on the circumference. Given an incident particle, denote by $x$ and $x^+$ the two points of intersection of the unfolded trajectory with the circumference. We are given $\measuredangle AOB = \ve$; therefore $x \in [0,\, \ve]$.

Denote by $\phi$ the angle between the direction vector of the unfolded trajectory and the radius at the first point of intersection; then the angle at the second point of intersection will be $-\phi$. Both angles are measured counterclockwise from the corresponding radius to the velocity; so, for example, $\phi > 0$ in figure \ref{fig unfold triangle}.

One has $x^+ = x + \pi - 2\phi$. The number of intersections of the unfolded trajectory with the images of the radii $OA$ and $OB$ coincides with the number of reflections of the {\it true} billiard trajectory and is equal to $n = n_\ve(x,\phi) = \lfloor \frac{x + \pi - 2\phi}{\ve} \rfloor$. In figure \ref{fig unfold triangle}, $n = 3$.

Denote by $\vphi$ and $\vphi^+$, respectively, the angles formed by the velocity of the {true} billiard trajectory with the outer normal to $AB$ at the moments of the first and second intersection with the opening $AB$. One easily sees that
\beq\label{app ineq}
|\vphi - \phi| \le \ve/2 \quad \text{and} \quad |\vphi^+ - (-1)^{n+1} \phi| \le \ve/2.
\eeq

The mapping $(x,\phi) \mapsto (\vphi,\vphi^+)$ defines a measure preserving one-to-one correspondence between a subspace of the space $[0,\, \ve] \times \interval$ with the measure $\frac{1}{2\sin(\ve/2)} dx \cdot \frac 12 \cos\phi\, d\phi$ and the space $\Box = \interval^2$ with the measure $\eta^\ve_\vee$. Consider also the mapping
$$
(x,\phi) \mapsto (\phi, (-1)^{n_\ve(x,\phi)+1} \phi)
$$
and the measure $\tilde\eta^\ve_\vee$ induced on $\Box$ by this mapping. One easily deduces from the inequalities (\ref{app ineq}) that the difference $\eta^\ve_\vee - \tilde\eta^\ve_\vee$ weakly converges to zero as $\ve \to 0$; therefore it is sufficient to prove the weak convergence
\beq\label{app convergence}
\tilde\eta^\ve_\vee \to \frac 12 (\eta_0 + \eta_\star) \quad \text{as} \ \ \ve \to 0.
\eeq

Introduce the function
$$
g_\ve(x,\phi) = \left\{ \begin{array}{rl} 1, & \text{if } \ n_\ve(x,\phi) \ \text{ is odd}\\
-1, & \text{if } \ n_\ve(x,\phi) \ \text{ is even} \end{array} \right..
$$
Similarly to the previous subsection \ref{app 1}, it suffices to prove that for any $-\pi/2 < \Phi_1 < \Phi_2 < \pi/2$,
\beq\label{app2 conv}
\lim_{\ve\to 0} \frac 1\ve \int_0^\ve\!\!\!\int_{\Phi_1}^{\Phi_2} g_\ve(x,\phi)\, \cos\phi\, d\phi\, dx = 0.
\eeq
Fix $x \in [0,\, \ve]$ and put $\phi_m = \frac 12 (x + \pi - m\ve)$. One has $g_\ve(x,\phi) = 1$, if $\phi_{2n-1} < \phi < \phi_{2n}$ and $g_\ve(x,\phi) = -1$, if $\phi_{2n} < \phi < \phi_{2n+1}$. We easily get that the integral $\int_{\Phi_1}^{\Phi_2} g_\ve(x,\phi)\, \cos\phi\, d\phi$ uniformly converges to zero as $\ve \to 0$ (actually, it is less than $2\ve$), and therefore, the convergence in (\ref{app2 conv}) also takes place.

\section*{Acknowledgements}

This work was partially supported by the Center for Research and Development in Mathematics and Applications (CIDMA) from ''Funda\c{с}\~ao para a Ci\^encia e a Tecnologia''{} FCT and by European Community Fund FEDER/POCTI, as well as by FCT: research project PTDC/MAT/72840/2006.

\end{document}